\documentclass[12pt]{article}
\title{{\bf An algorithm for de Rham cohomology groups
            of the complement of an affine variety 
            via $D$-module computation}}
\author{Toshinori Oaku and Nobuki Takayama}
\date{January 26, 1998}

\def\Ann{ {\rm Ann}\, }
\def\pd#1{ \partial_{#1} }
\def\C{{\bf C}}
\def\hol{ {\cal O} }
\def\Hom{ {\cal H}{\rm om}\, }
\def\N{ {\bf N} }
\def\Q{{\bf Q}}
\def\R{ {\bf R}}
\def\V{ V }
\def\Z{ {\bf Z} }

\def\Dsc{ {\cal D} }
\def\Hsc{ {\cal H} }
\def\Isc{ {\cal I} }
\def\Ksc{ {\cal K} }
\def\Lsc{ {\cal L} }
\def\Msc{ {\cal M} }
\def\Nsc{ {\cal N} }
\def\Osc{ {\cal O} }
\def\Psc{ {\cal P} }
\def\Vsc{ {\cal V} }

\newenvironment{proof}{\par\noindent Proof:\ }{${\tt [}\kern-0.2mm{\tt ]}$}

\newtheorem{theorem}{Theorem}[section]
\newtheorem{algorithm}[theorem]{Algorithm}
\newtheorem{corollary}[theorem]{Corollary}
\newtheorem{definition}[theorem]{Definition}
\newtheorem{example}{Example}[section]
\newtheorem{lemma}[theorem]{Lemma}
\newtheorem{procedure}[theorem]{Procedure}
\newtheorem{proposition}[theorem]{Proposition}
\newtheorem{remark}{Remark}[section]

\def\comment#1{  }

\begin{document}
\maketitle

\comment{ Abstract :
We give an algorithm to compute the following cohomology groups
on $U = \C^n \setminus V(f)$ for any non-zero polynomial
$f \in \Q[x_1, \ldots, x_n]$;
1. $H^k(U, \C_U)$, $\C_U$ is the constant sheaf on $U$ with stalk $\C$.
2. $H^k(U, \Vsc)$, $\Vsc$ is a locally constant sheaf of rank $1$ on $U$.
We also give partial results on computation of cohomology groups on $U$
for a locally constant sheaf of general rank and on computation of
$H^k(\C^n \setminus Z, \C)$ where $Z$ is a general algebraic set.
Our algorithm is based on computations of Gr\"obner bases in the ring
of differential operators with polynomial coefficients.
14F40, 14Q99, 55N30.
}

\setcounter{section}{-1}
\section{Introduction}

In this paper, 
we give an algorithm to compute the following cohomology groups
on $U = \C^n \setminus V(f)$
for any non-zero polynomial
$f \in \Q[x_1, \ldots, x_n]$;
\begin{enumerate}
\item $H^k(U, \C_U)$, $\C_U$ is the constant sheaf on $U$ with stalk $\C$.
\item $H^k(U, \Vsc)$, $\Vsc$ is a locally constant sheaf of rank $1$ on $U$.
\end{enumerate}
We also give partial results on computation of cohomology groups on $U$
for a locally constant sheaf of general rank and on computation of
$H^k(\C^n \setminus Z, \C)$ where $Z$ is a general algebraic set.

Our algorithm is based on computations of Gr\"obner bases in the ring
of differential operators with polynomial coefficients,
algorithms for functors in the theory of $\Dsc$-modules 
(\cite{OakuAdvance} and \cite{OakuAdvance2}),
and Grothendieck-Deligne comparison theorem \cite{G}, \cite{D},
which relates sheaf cohomology groups and de Rham cohomology groups.

One advantage of the use of the ring of differential operators
in algebraic geometry is that, for example, 
$\Q[x,1/x]$, which is the localized module of $\Q[x]$ along $x$,
is not finitely generated as a $\Q[x]$-module, but it can be regarded
as a finitely generated $\Q \langle x, \pd{x} \rangle$-module 
where $\pd{x} = \partial/\partial x$.
In fact, we have
$ 1/x^k = (-1)^{k-1}(1/(k-1)!) 
  \left( {\partial \over {\partial x}} \right)^{k-1} {1 \over x}$.
Computation of the localization of a given $\Dsc$-module and computation
of the integration functor are the most important part of our algorithm.

\medbreak

As an introduction to this paper,
it will be the best to mention about
how we started this project.

Let us consider the differential equation for the function
$f = (x-u)^a(x-v)^b $
where $u$ and $v$ are rational numbers ($u < v$) and $a$ and $b$ are
rational numbers such that $a+b \not\in \Z$.
The function $f$ satisfies the differential equation
$$ p f = 0,  \  p = (x-u) (x-v) \pd{x} - a (x-v) - b(x-u). $$
Let ${\hat p}$ be the formal Fourier transform of this operator
\begin{eqnarray*}
 {\hat p} &=& (-\pd{x}-u) (-\pd{x}-v) x - a (-\pd{x}-v) - b (-\pd{x}-u)  \\
          &=& x \pd{x}^2 + (ux+vx+2+a+b)\pd{x} + uvx + u + v + av + bu 
\end{eqnarray*}
and $A$ be the ring of differential operators $\Q \langle x, \pd{x} \rangle$.
We want to evaluate the dimension of the $\Q$-vector space
$$ A/(A {\hat p} + x A) \simeq (A/A {\hat p})/x (A/A {\hat p}) \simeq
  (A/A {\hat p})/ {\rm Im}\, x \cdot, $$
which is called  {\it the {\rm (}$0$-th{\rm )} restriction} of $A/A {\hat p}$ along
$x=0$.
Now, we can apply the algorithm for the $\Dsc$-module theoretic restriction 
in \cite[Section 5]{OakuAdvance}
to evaluate the dimension.  Here, we need what is called a $b$-function 
for the evaluation,  
and it is nothing but the indicial (characteristic) polynomial at $x=0$ of
the ordinary differential operator ${\hat p}$ that appears in the classical
method of Frobenius; 
the polynomial is $ s (s - a - b)$.
Applying Proposition 5.2 in \cite{OakuAdvance} with this $b$-function,
we conclude that the dimension is equal to one,
which is the number of the bounded segments of
$ \R \setminus \{u, v\}$.

Next, we tried to evaluate the dimension of
$ A_2 / (A_2 {\hat p} + A_2 {\hat q} + x A_2 + y A_2) $
where $A_2$ is the ring of differential operators generated by
$x, y, \pd{x}$ and $\pd{y}$, 
and
$p$ and $q$ are differential operators which annihilate the function
$f = x^a y^b (1-x-y)^c$;
$ p = x(1-x-y) (\pd{x} - a/x + c/(1-x-y)) $ and
$ q = y(1-x-y) (\pd{y} - b/y + c/(1-x-y)) $.
We evaluate the dimension, this time with a computer program \cite{Kan-sm1},
by iterating to apply the algorithm 
for computing the $0$-th restriction in \cite{OakuAdvance} 
firstly for $x$ and secondly for $y$.
The result is again one, which is equal to the number of
the bounded cells of the hyperplane arrangement
$ \R^2 \setminus \{ (x,y)\,|\, x y (1-x-y) = 0  \}$.
It is well known in the theory of hypergeometric functions that
the number of bounded cells is equal to the dimension of
the middle dimensional twisted cohomology group associated with $f$, 
which is equal to
the rank of the corresponding hypergeometric system.
(Strictly speaking, it turns out that $p$ and $q$ do not generate 
the annihilating ideal $\{\ell\in A_2 \mid \ell f = 0\}$ for $f$
(see Example \ref{example:example4.1}); 
however, the ideal generated by $p$, $q$ 
happens to be `close enough' to the annihilating ideal.)

Inspired by the observation above, we started the project 
to obtain an algorithm for computing the cohomology groups
of the complement 
of an affine variety by elaborating the method sketched above.

\section{Computation of cohomology groups of the complement of an affine hypersurface}

For any non-zero polynomial $f(x) \in \Q[x_1, \ldots, x_n]$,
we  prove the following theorem.

\begin{theorem} \label{theorem:constantSheaf}
Put $X = \C^n$, $Y = V(f) := \{x \in X \mid f(x)= 0\}$, 
and $U = X \setminus Y$.   
Then the cohomology group $H^k(U,\C_U)$ is computable for any integer $k$, 
where $\C_U$ denotes the constant sheaf on $U$ with stalk $\C$. 
\end{theorem}

Note that  $H^k(U,\C_U)$ is the $k$-th cohomology group (with coefficients in 
$\C$) of the $2n$-dimensional real $C^\infty$-manifold $U_{cl}$ underlying
$U$.
Also note that $H^k(U,\C_U) = 0$ for $k > n$ since $U$ is affine
and that $H^0(U,\C_U) = \C$ since $U$ is connected.

In the theorem above, we may replace $\Q$ by any computable field.
Here, we mean by a computable field a subfield $K$ of $\C$ such that 
each element of $K$ can be expressed by a finite set of data so that
we can decide whether two such expressions correspond to the same element, 
and that the addition, subtraction, multiplication, and division 
in $K$ are computable by the Turing machine.
For example, any algebraic extension field of $\Q$ of finite rank is
a computable field by virtue of Gr\"obner bases and factorization 
algorithms over algebraic number fields. 

In this section, we illustrate an algorithm to compute 
the cohomology groups.
Correctness will be proved as a special case of the corresponding theorem
for cohomology groups with the coefficients in a locally constant
sheaf of rank one.
In order to compute the cohomology groups,
we translate the problem to that of
computations of functors, especially to that of the de Rham functor,
of $A_n$-modules,
which are studied in a series of papers \cite{OakuAdvance}, \cite{OakuAdvance2}.
Here, $A_n$ is the ring of differential operators with polynomial
coefficients and is called {\it the Weyl algebra}.
The computations of functors are based on the Buchberger algorithm
to compute Gr\"obner bases in the Weyl algebra.
We shall quickly review the definition of Weyl algebra and
the Gr\"obner basis.
See \cite{Galligo}, \cite{Castro}, \cite{TakayamaJJAM}, \cite{OakuBook} for 
details, \cite{Sturmfels-Takayama} for an introduction,
and \cite{Kan-sm1} for implementations.

The Weyl algebra
$$A_n = \Q \langle x_1, \ldots, x_n, \pd{1}, \ldots, \pd{n} \rangle  $$ 
is the ring of non-commutative polynomials
generated by $2n$ elements $x_i, \pd{i}$, 
$(i=1, \ldots, n)$
satisfying the relations
\begin{eqnarray*}
 x_i x_j &=& x_j x_i, \ \pd{i} \pd{j} = \pd{j} \pd{i} \\
 \pd{i} x_j &=& x_j \pd{i} + \left\{
 \matrix{ 1 &  (i = j) \cr
          0 &  (i \not= j). \cr
 }
 \right. \\
\end{eqnarray*}

The theory (and practice) of Gr\"obner bases works perfectly well
for left ideals in the Weyl algebra $A_n$. We quickly review the
relevant basics. Every element $p$ in $A_n$ can be written uniquely
as a $\Q$-linear combination of 
{\it normally ordered} monomials
$x^a \pd{}^b$.
This representation of $p$ is called the {\it normally ordered representation}.
For example, the monomial $\, \pd{1} x_1 \pd{1} \,$ is not normally ordered.
Its normally ordered representation is  $\, x_1 \pd{1}^2 + \pd{1}$.

Consider the commutative polynomial ring in $2n$ variables
$$ {\rm gr\,}(A_n) \quad =\quad \Q [x_1, \ldots , x_n, \xi_1,  
\ldots , \xi_n] $$
and the  $\Q$-linear map
$\, {\rm gr\,} : A_n \rightarrow {\rm gr\,}(A_n) ,\, 
  x^a \pd{}^b \mapsto  x^a \xi^b $.
Let $<$ be any term order on ${\rm gr\,}(A_n)$.
This gives  a total order among normally ordered monomials in $A_n$ via
$\,  x^A \pd{}^B > x^a \pd{}^b \ \Leftrightarrow
    x^A \xi^B > x^a \xi^b $.
For any element $p \in A_n$ let  $in_<(p)$ denote 
the highest monomial $ x^A \pd{}^B$ in the normally 
ordered representation of $p$. 
If $I$ is a left ideal in $A_n$ then its
{\it initial ideal} is the ideal ${\rm gr\,}(in_<(I)) $ in  ${\rm gr\,}(A_n)$
generated by all monomials $\,{\rm gr\,}(in_<(p))\,$ 
for $p \in I$.  Clearly,  ${\rm gr\,}(in_<(I))$ is generated by finitely many 
monomials $  x^a \xi^b $. 
A finite subset $G$ of $I$ is called {\it a Gr\"obner basis } of $I$
with respect to the term order $<$ if 
$\{ {\rm gr}\, (in_< (q)) \,|\, q \in G \}$
generates
${\rm gr}\,(in_< (I))$.
Noting that
 $ in_<(p) \leq in_<(q)$ implies
$ in_<(hp) \leq in_<(hq) $  for all $h \in A_n$,
one proves that the reduced Gr\"obner basis of $I$ is unique and finite, 
and can be computed using Buchberger's algorithm.
Any left (or right) ideal in $A_n$ is finitely generated and we denote by
$ \langle p_1, \ldots, p_m \rangle $
the left ideal in $A_n$ generated by
$p_1, \ldots, p_m \in A_n$.

Most constructions in the commutative algebra can be reduced to
computations of Gr\"obner bases.
This is also the case with some constructions for modules 
over the Weyl algebra.
For example,
the construction of a free resolution of a left coherent $A_n$-module 
is a straightforward generalization of algorithms of constructing 
free resolutions of modules over the ring of polynomials.
As to algorithms to construct a free resolution by {\it the Schreyer order},
see \cite[p.167 Theorem 3.7.13]{AdamsBook}, \cite[Theorem 15.10]{EisenbudBook}
and \cite{StillmanResolution}.
We note that 
 an algorithm to compute a sheaf
cohomology on the $n$-dimensional projective space is given
based on computation of syzygies in the ring of polynomials
(\cite{EisenbudCohomology}).
Computation of {\it an elimination ideal} in the Weyl algebra is also a 
straightforward generalization of computation of an elimination 
ideal in the ring of polynomials,
(see, e.g., \cite[p.69 Theorem 2.3.4]{AdamsBook},
\cite[p.114 Theorem 2]{CoxBook}).
These two constructions will be used in our algorithm to obtain
cohomology groups (see Algorithm \ref{algorithm:cohomology-top} Step 3,
Procedure \ref{procedure:annihilator} Step 2 
and Procedure \ref{procedure:severalAnnihilator} Step 2).

However, the non-commutativity 
causes some difficulty in constructing various objects in the category
of modules over the Weyl algebra.
For example, to compute
the tensor product of right and left $A_n$-modules in the derived category,
special care must be taken. 
This problem has been an open problem since \cite{TakayamaISSAC90}. 
As a special (but important) case of the tensor product computation as above, 
we give 
an algorithm for the $\Dsc$-module theoretic restriction of an $A_n$-module
by using the $V$-filtration and the $b$-function (or the indicial 
polynomial).  
As to details, see 
\cite[Proposition 5.2, Theorem 5.7, Algorithm 5.10]{OakuAdvance} 
and \cite{OakuAdvance2}.
U.Walther \cite{Walther} solved a related problem of computing 
algebraic local cohomology groups based on $V$-filtration, $b$-function
and the Cech complex.
As we will see in Section \ref{section:acyclic},
his algorithm gives an algorithm for 
Theorem \ref{theorem:constantSheaf} different from ours explained below.

We have explained a general background on an algorithmic
treatment of modules over the Weyl algebra.
Now let us  explain our algorithm to compute cohomology 
groups by a top-down expansion.

\begin{algorithm}
\rm 
(Computation of the cohomology groups $H^k(U,\C_U)$)
\label{algorithm:cohomology-top}

\noindent
Input : a polynomial $f \in \Q[x_1, \ldots, x_n]$.\\
Output : $H^k(U,\C_U)$ for $0 \leq k \leq n$ where $U = \C^n \setminus V(f)$.

\begin{enumerate}
\item  \label{algorithm:cohomology-top:step1}
Find a left ideal $I$ such that 
$$  \Q [x, {1 \over f}] \simeq A_n/I $$
as a left $A_n$-module.
\item \label{algorithm:cohomology-top:step2}
Let $J$ be the formal Fourier transform of $I$;
$$  J = I_{|_{ x_i \mapsto -\pd{i}, \pd{i} \mapsto x_i}}. $$
\item Compute a free resolution of length $n+1$ 
$$    A_n^{p_{-(n+1)}}
      \stackrel{\cdot L^{-(n+1)}}{\longrightarrow}    A_n^{p_{-n}}
      \stackrel{\cdot L^{-n}}{\longrightarrow} A_n^{p_{-(n-1)}}
      \cdots 
      \stackrel{\cdot L^{-1}}{\rightarrow} 
      A_n^{p_0} \rightarrow A_n/J \rightarrow 0,
$$
($p_0 = 1$)
of $A_n/J$ by using Schreyer's theorem \cite[Theorem 15.10]{EisenbudBook}
with an order which refines the partial order defined by the weight vector
$$
\matrix{
& \partial_1 & \cdots &   \partial_{n} & x_1 & \cdots &  x_n &\cr
(&  1     & \cdots  &       1           & -1  & \cdots &  -1 &). \cr
}
$$
\item Compute the cohomology groups of the complex of $\Q$-vector
spaces
$$    ( A_n/(x_1 A_n + \cdots + x_n A_n) \otimes_{A_n}
        A_n^{p_{-k}} ,  \stackrel{1 \otimes L^{-k}}{\longrightarrow}).
$$
Then, the $(k-n)$-th cohomology group 
$ {\rm Ker}\,(1 \otimes L^{k-n})/ {\rm Im}\,(1 \otimes L^{k-n-1})$
of the complex above tensored with
$\C$ gives $H^k(U,\C_U)$.
\end{enumerate}
\end{algorithm}

The step 1 will be explained in Procedure \ref{procedure:annihilator}
in detail and the steps 2, 3 and 4 will be explained in
Procedure \ref{procedure:integral} in detail.

We note that the steps 2, 3 and 4 are nothing but the computation of
$$ H^{k-n}(A_n/(\pd{1} A_n + \cdots + \pd{n} A_n) \otimes_{A_n}^{L}
           A_n/I),
$$
which is denoted by 
$$ \int^{k-n}_{\C^n} A_n/I   = \int_{\C^n}^{k-n} \Q[x,1/f]
$$ 
in the theory of ${\cal D}$-modules.
We shall prove that this cohomology group tensored with $\C$ is equal to
$H^k(U,\C_U)$ by the Grothendieck-Deligne comparison theorem
in Section \ref{deRhamCohomologyAndIntegrationOfDmodules}.
Here, for a left $A_n$-module $M$ and a right $A_n$-module $N$, 
we denote by $N \otimes_{A_n}^{L} M$  the complex
$$ ( N \otimes_{A_n} M^i , 1 \otimes d^{i-1} ; i=1, 0, -1, -2, \ldots ) $$
where
$$   \cdots \stackrel{d^{i-1}}{\longrightarrow} M^i
            \stackrel{d^{i}}{\longrightarrow}  \cdots
            \stackrel{d^{-2}}{\longrightarrow}  M^{-1}
            \stackrel{d^{-1}}{\longrightarrow}  M^{0}
            \longrightarrow M \longrightarrow 0, \  
   \mbox{ (exact) }
$$
$M^i$ is a free $A_n$-module, $d^i$ is a left $A_n$-morphism
and
$M^1 = 0$, $d^0 = 0$.
It is known that there exists a finite length free resolution
for a given finitely generated left $A_n$-module $M$
(e.g., apply the method of \cite[p.336 Corollary 15.11]{EisenbudBook}
to our case).

\begin{remark} \rm
For a left $A_n$-module $M=A_n/I$,
the left $A_{n-1}$-module
$(A_n/\pd{n} A_n) \otimes_{A_n} M = A_n/(I+\pd{n}A_n)$
is called the $0$-th integral of $M$ with respect to $x_n$.
Why is it called the integral?
Let us explain an intuitive meaning of this terminology.

Let $f$ be a function of $x_1, \ldots, x_n$.
We suppose that the function $f$ is rapidly decreasing with respect
to the variable $x_n$ and put
$I = \Ann f = \{ \ell \in A_n \,|\, \ell f = 0 \}$.
Then the $A_n$-module generated by the function $f$ is isomorphic
to the left $A_n$-module $A_n/I$.
Put
$ g(x_1, \ldots, x_{n-1}) = \int_{-\infty}^\infty f(x_1, \ldots, x_n) dx_n$.
Then, we have
$ [(I + \pd{n} A_n) \cap A_{n-1}] g = 0 $.
In fact, since any element $\ell$ in $(I+\pd{n}A_n) \cap A_{n-1}$
can be written as
$ \ell = \ell_1 + \pd{n} \ell_2$,
$ \ell_1 \in I,  \ell_2 \in A_n$,
we have
$ \ell g = \int_{-\infty}^\infty (\ell_1 + \pd{n} \ell_2) f dx_n
 = \int_{-\infty}^\infty \pd{n}(\ell_2 f) dx_n = 0$.
Therefore, $g$ can be regarded as a solution of the differential equations 
corresponding to the left $A_{n-1}$-submodule 
$A_{n-1}/(I+\pd{n}A_n)\cap A_{n-1}$ of
$A_n/(\pd{n} A_n) \otimes_{A_n} M$.
Note that $A_n/(\pd{n} A_n) \otimes_{A_n} M$ itself describes a system of 
differential equations for 
$\int_{-\infty}^\infty x_n^j f\,dx_n$ with $j \geq 0$.  
\end{remark}

Let us explain in detail the step 1 of Algorithm 
\ref{algorithm:cohomology-top}.
This algorithm is given in \cite{OakuMEGA96}

\begin{procedure} \rm \label{procedure:annihilator}
(Computing the differential equations for $1/f^{-r_0}$;
step 1 of Algorithm \ref{algorithm:cohomology-top}). \\

\noindent
Input: $f$. \\
Output: a left ideal $I$ of $A_n$ such that $\Q[x,1/f] \simeq A_n/I$.

\begin{enumerate}
\item (Computation of the annihilating ideal of $f^s$) 

Compute
$$
\langle t-f(x), {{\partial f} \over
                 {\partial x_1}} \pd{t} + \pd{1}, \ldots,
                {{\partial f} \over
                 {\partial x_n}} \pd{t} + \pd{n} \rangle \cap
\Q[ t \pd{t}]\langle x, \pd{x}\rangle.
$$
Replacing $t \pd{t}$ by $-s-1$, we obtain the left ideal
$ \Ann f^s  $ in $\Q[s]\langle x, \pd{x} \rangle$.
(Call Procedure \ref{procedure:intersectionWithSubring} with $d=1$ 
to compute the intersection
of the left ideal and the subring $\Q[t \pd{t}] \langle x, \pd{x}\rangle$.)

\item (Computation of the $b$-function of $f$)

Compute the generator $b(s)$ of
$$  \langle \Ann f^s, f \rangle \cap \Q[s] $$
by an elimination order $ x, \pd{x} > s $.
\item
Let $r_0$ be the minimum integral root of $b(s)=0$.
Put $I = (\Ann f^s)_{_{ s \rightarrow r_0}}$.
Then, we have $\Q[x,{1 \over f}] \simeq A_n/I$.
\end{enumerate}
\end{procedure}

The polynomial $b(s)$ is called the (global) 
{\it Bernstein-Sato polynomial} or the {\it $b$-function} of $f$, which 
coincides with the (global) {\it indicial polynomial} of $\delta(t-f)$,
and plays an important role in our algorithms.
This polynomial is the minimal degree polynomial satisfying the relation
$$  L f^{s+1} = b(s) f^s, \quad \exists L \in \Q[s]\langle x, \pd{x} \rangle.
$$
The left module $A_n/I$ is a {\it holonomic} $A_n$-module 
(or called  a module belonging to the Bernstein class).
The holonomicity of $A_n/I$ and the existence of the $b$-function
were shown by I.N.Bernstein.
See, e.g., \cite{Bernstein} and \cite[p.13, 5.5 Theorem]{BjorkBook}.
It is known that when $f \not= {\rm const}$, $b(s)$ always 
has $s+1$ as a factor.
M.Kashiwara proved that
all the roots of $b(s) = 0$ are negative rational numbers
for any $f \in \C[x_1, \ldots, x_n]$ \cite{KashiwaraBfunction}.

\begin{remark} \rm
The left $A_n$-isomorphism
$  \Q [x,1/f] \stackrel{\varphi}{\longrightarrow} A_n/I $
is expressed as $\varphi = \varphi_2 \circ \varphi_1$ where
the left $A_n$-isomorphisms $\varphi_1$ and $\varphi_2$
$$ \Q [x,1/f] \stackrel{\varphi_1}{\longrightarrow} A_n f^{r_0}
              \stackrel{\varphi_2}{\longrightarrow} A_n/I
$$
are defined as
$$ \varphi_2(f^{r_0}) = 1 \in A_n/I, \quad
   \varphi_1(f^{r_0}) = f^{r_0} \in A_n f^{r_0}
$$
and
$$  \varphi_1(f^{r_0-k}) =
 { {L(r_0-k) \cdots L(r_0-1)} \over {b(r_0-k) \cdots b(r_0-1)}} f^{r_0},
\  (k=1, 2, \ldots ).
$$
Hence, for example, we have
$$\varphi^{-1}(\pd{i}) = r_0 (\partial f/\partial x_i) f^{r_0-1}
= \varphi_1^{-1}(r_0 (\partial f/\partial x_i) L(r_0-1) f^{r_0}/b(r_0-1)).
$$
\end{remark}

\begin{example}  \rm
For $f = x (1-x)$, we have
$$ \Ann f^s = \langle x (1-x) \pd{x} - s (1-2x) \rangle. $$
The $b$-function of $f$ is $s+1$ with
$$ ((1-2x) \pd{x} + 4 (1+s)) f^{s+1} = (s+1) f^s, $$
and hence we get 
$$ \Q[x,1/f] \simeq \Q\langle x, \pd{x} \rangle 
            / \langle x(1-x) \pd{x} + (1-2x) \rangle.
$$
\end{example}

\begin{example} \rm
Put $f = x^3-y^2$.  We compute the left ideal $I$ such that
$\Q[x,y,1/f] \simeq A_2/I$.
Here is a log of the output of {\tt kan/k0}, which may be self-explanatory.
The system {\tt k0} is a translator 
that compiles Java like inputs to codes for 
{\tt kan/sm1}, which is a Postscript like language for computations in
the ring of differential operators \cite{Kan-sm1}.
\begin{verbatim}
In(9)= a = annfs(x^3-y^2,[x,y]);
Computing the Groebner basis of 
[ v*t+x^3-y^2 , -v*u+1 , -3*u*x^2*Dt+Dx , 2*u*y*Dt+Dy ] 
with the order u, v > other elements.

In(10)=a :
[3*x^2*Dy+2*y*Dx , -6*(-1-s)-2*x*Dx-3*y*Dy-6] 

In(11)=b=ReducedBase(Eliminatev(Groebner(Append[a,y^2-x^3]),
                                [x,y,Dx,Dy]));
In(12)= b:
[   -216*s^3-648*s^2-642*s-210 ]
In(13)=Factor(b[0]):
[ [ -6 , 1 ], [6*s+5 , 1 ], [6*s+7 , 1 ], [s+1 , 1 ]  ] 
\end{verbatim}
Since $s = -1$ is the minimum integral root of the $b$-function, we have
$$ \Q[x,y,1/f] \simeq A_2/\langle 3x^2\pd{y}+2y\pd{x}, 
                            -2x\pd{x}-3y\pd{y}-6 \rangle.
$$
\end{example}

Finally, let us explain our algorithm for computing 
$$ H^{k-n}(A_n/(\pd{1} A_n + \cdots + \pd{n} A_n) \otimes_{A_n}^{L}
           A_n/I).
$$
This is a detailed explanation of steps 2, 3 and 4 of Algorithm
\ref{algorithm:cohomology-top}.
We can compute the cohomology groups by applying 
\cite[Theorem 5.3]{OakuAdvance2} to the Fourier transformed ideal
$J$ of $I$.
Correctness will be discussed in Sections \ref{deRhamCohomologyAndIntegrationOfDmodules} 
and \ref{section:integration}.

Put
$$ 
\matrix{
&&   & \partial_1 & \cdots &   \partial_{n} & x_1 & \cdots &  x_n &\cr
w&=&(&  1     & \cdots  &       1           & -1  & \cdots &  -1 &). \cr
}
$$
and
$$ F_k = \{ f \in A_n \,|\, ord_w(f) \leq k \} $$
where
$$ ord_w(x^a \pd{}^b) = -|a|+|b|. $$
$\{F_k\}$ is called the $V$-filtration.

\begin{procedure}
\rm \label{procedure:integral}
\cite{OakuAdvance2} (Computing the $D$-module theoretic integral of $A_n/I$;
steps 2, 3 and 4 in Algorithm \ref{algorithm:cohomology-top})

\noindent
Input: a left ideal $I$ of $A_n$. ($A/I$ is holonomic.)\\
Output: The $-k$-th cohomology groups of 
$ A_n/(\pd{1} A_n + \cdots + \pd{n} A_n) \otimes_{A_n}^L A_n/I $ 
for $0 \leq k \leq n$.

\begin{enumerate}
\item
Let $J$ be the formal Fourier transform of $I$;
$$  J = I_{|_{ x_i \mapsto -\pd{i}, \pd{i} \mapsto x_i, (i=1, \ldots, n)}}. $$

\item  Let $G$ be a Gr\"obner basis of the left ideal $J$
with the weight vector $w$.
Find the generator
$ b(\theta_1 + \cdots + \theta_n) $
of
$$\langle in_w(G) \rangle \cap \Q[\theta_1+ \cdots + \theta_n], \quad
  \theta_i = x_i \pd{i}.
$$
\item
Let $k_1$ be the maximum integral root of $b(s)=0$.
If there exists no integral root, then quit; 
the cohomology groups are all zero in that case. 
\item
Let $<_w$ be a refinement of the partial order by $w$.
Construct a free resolution 
$$    A_n^{p_{-(n+1)}}
      \stackrel{\cdot L^{-(n+1)}}{\longrightarrow}    A_n^{p_{-n}}
      \stackrel{\cdot L^{-n}}{\longrightarrow} A_n^{p_{-(n-1)}}
      \cdots 
      \stackrel{\cdot L^{-1}}{\rightarrow} 
      A_n^{p_0} \rightarrow A_n/J \rightarrow 0
$$
with $p_0 = 0$ 
by using the Schreyer orders associated with $<_w$.  
\item (Computation of degree shifts)
Put $s_1^0 = 0$ and
$$
  s_i^{k+1} = \max_{1\leq j\leq p_{-k}}\left(
      ord_w(L^{-(k+1)}_{ij}) + s_j^k  \right) \quad (1 \leq i \leq p_{-(k+1)})
$$
successively. 
\item
Compute the cohomology groups of the induced complex 
$$
\begin{array}{ccc}
 \cdots 
 &\stackrel{\cdot {\bar L}^{-2}}{\longrightarrow}&
 F_{k_1-s_1^1}/( F_{-1} + x A_n ) \bigoplus \cdots \bigoplus
 F_{k_1-s_{p_{-1}}^1}/( F_{-1} + x A_n)          \\
&\stackrel{\cdot {\bar L}^{-1}}{\longrightarrow}&
 F_{k_1}/ ( F_{-1} + xA_n )
 \stackrel{\bar L^0}{\rightarrow} 0 \\
\end{array}
$$
as a complex of $\Q$-vector space
where $x A_n = x_1 A_n + \cdots + x_n A_n $.
Then, the $(k-n)$-th cohomology group
$$ {\rm Ker}\, {\bar L}^{k-n}/
   {\rm Im}\, {\bar L}^{k-n-1}$$
of this complex tensored by ${\bf C}$ gives 
$H^k(U,\C_U)$.
\end{enumerate}
\end{procedure}

In step 2, we denote by
$ in_w(\sum_{(a,b) \in I} c_{ab} x^a \pd{}^b ) $
the $w$-leading form 
$$ \sum_{ \langle w, (a,b) \rangle = m} c_{ab} x^a \pd{}^b, \quad
   m = \max_{ (a,b) \in I} \langle w, (a,b) \rangle
$$
where $\langle \cdot, \cdot \rangle$ denotes the standard inner product
in $\Z^{2n}$.
Put $in_w(G) = \{ in_w(g) \,|\, g \in G \}$.
One needs to compute the intersection of the left ideal
$in_w(J) = \langle in_w(G) \rangle $
and the subring.  This can be done in a procedure similar to the one explained
in Section \ref{section:intersection}.

In step 4, we compute Gr\"obner bases with the Schreyer orders over
the order $<_w$ to construct a resolution.
Note that $<_w$ is not a well-order, which causes a difficulty of computation
of Gr\"obner basis in our non-commutative situation.
There are two ways to overcome this difficulty;
one is to use the $F$-homogenization introduced in \cite{OakuJMS} (see also
\cite[Section 3]{OakuBfunction}) and the other is the use of
the homogenized Weyl algebra which has a homogenization variable $h$ so that
the relation $\pd{i} x_i = x_i \pd{i}+h^2$ holds.
The homogenized Weyl algebra was introduced in {\tt kan/sm1} 
\cite{Kan-sm1} since version 2 released in 1994.
See \cite{ACG} on a theoretical study on this homogenization technique.

In Step 6, we truncate the complex from above by using $k_1$; 
we could also make truncation from below by using the minimum integral 
root of $b(s)=0$, which would somewhat reduce the complexity of Step 6.  
See \cite{OakuAdvance2} for details.
Moreover, if we need to compute $H^k(U,\C_U)$ only for $k \geq \ell$, 
then a resolusion of length $n-\ell+1$ suffices in Step 4.

\begin{example} \rm
We take $f = x(1-x)$ in $\Q[x]$. We denote by $A$ the Weyl algebra
$\Q\langle x, \pd{x} \rangle$.
As we have seen, we have
$ \Q[x,1/f] \simeq A/\langle p \rangle , \quad p = x (1-x) \pd{x} - (2x-1)$.
The formal Fourier transform of $p$ is
$ {\hat p}= -x\pd{x}^2-x \pd{x}$.
By multiplying $-x$ from the left, we have
$$ -x{\hat p} = \theta(\theta-1) + x \theta , \quad \theta = x \pd{x}.$$
Therefore the $b$-function is equal to $s(s-1)$ and
$k_1 = 1$.
The resolution of $A/\langle {\hat p}\rangle$ is
$$
 0 \longrightarrow A
  \stackrel{\cdot (x \pd{x}^2-x \pd{x})}{\longrightarrow}
  A
  \longrightarrow
  A/\langle {\hat p}\rangle 
  \longrightarrow 0.
$$
Since $k_1=1$ and  the degree shift by
$x \pd{x}^2-x \pd{x}$ is equal to $1$,
the truncated complex is 
$$
 0 \longrightarrow
 F_0/( F_{-1} + xA )
 \stackrel{\cdot (x \pd{x}^2-x \pd{x})}{\longrightarrow}
 F_1/( F_{-1} + xA )
 \longrightarrow 0.
$$
Since
$$ F_0/( F_{-1}+ xA ) = \Q, \ 
   F_1/( F_{-1} + xA ) = \Q + \Q \pd{x}
$$
and
$ 1 \cdot (x \pd{x}^2-x \pd{x}) \equiv 0 $
in 
$  F_1/( F_{-1} + xA ) $,
we conclude that
$$ H^{-1} = F_0 / ( F_{-1} + xA ) = \Q, 
 \quad  H^{0} = F_0 / ( F_{-1} + xA ) = \Q^2. $$
Hence, the cohomology groups of $U= \C \setminus \{0, 1\}$
are
$$ H^0(U,\C_U) = \C, \quad H^1(U,\C_U) = \C^2. $$
The two generators of $H^1$ correspond to two loops that encircle
the points $x=0$ and $x=1$ respectively in view of the Poincar\'e duality
of homology groups and cohomology groups.
\end{example}

\begin{example} \rm
(Cohomology groups of $\C^2 \setminus V(x^3-y^2)$)

This is the output of {\tt kan/k0}.
\begin{verbatim}
In(43)= bb=bfunctionForIntegral([3*x^2*Dy+2*y*Dx,-2*x*Dx-3*y*Dy-6],[x,y]);
In(44)= bb:
[    -216*s^3+432*s^2-264*s+48 ] 
In(45)= Factor(bb):
[[-24 , 1 ], [ 3*s-2 , 1 ], [ s-1 , 1 ], [ 3*s-1 , 1 ]] 
In(46)=integralOfModule([3*x^2*Dy+2*y*Dx , -2*x*Dx-3*y*Dy-6],[x,y],1,1,2):
\end{verbatim}
Here, $1,1,2$ specify the mimimum and the maximum integral roots, and 
the length of the resolution respectively.
\begin{verbatim}
0-th cohomology:  [    0 , [   ]  ] 
-1-th cohomology:  [    1 , [   ]  ] 
-2-th cohomology:  [    1 , [   ]  ] 
\end{verbatim}
The output means that
$$ H^0(U,\C_U) = \C, \quad H^1(U,\C_U) = \C, \quad H^2(U,\C_U) = 0.$$

Let us explain this example a little more precisely.
For 
$f = x^3-y^2$, we have
$\Q[x,y,1/f] \simeq A_2/I$ with
$$
I = \langle 2x\pd{x}+3y\pd{y}+6, \,\,  3x^2\pd{y}+2y\pd{x} \rangle.
$$
Its Fourier transform is $A_2/J$ with 
$$
J = \langle -2x\pd{x}-3y\pd{y}+1, \,\,  3y\pd{x}^2-2x\pd{y} \rangle.
$$
The $b$-function of $A_2/J$ is $(s-1)(3s-1)(3s-2)$.  
Hence we put $k_1 = 1$. 
The free resolution of $A_2/J$ is given by 
$$    
      0 \rightarrow A_2
      \stackrel{\cdot L^{-3}}{\longrightarrow} A_2^{4}
      \stackrel{\cdot L^{-2}}{\longrightarrow} A_2^{4}
      \stackrel{\cdot L^{-1}}{\rightarrow} 
      A_2 \rightarrow A_2/J \rightarrow 0
$$
with 
\begin{eqnarray*}
L^{-1} &=& \left( \matrix{
              &  -2x\pd{x}-3y\pd{y}+1 & \cr
              &   3y\pd{x}^2-2x\pd{y} & \cr
              &  -9y^2\pd{y}\pd{x}-3y\pd{x}-4x^2\pd{y} & \cr 
              &  -27y^3\pd{y}^2-27y^2\pd{y}+3y+8x^3\pd{y} & \cr
        }\right), \\
L^{-2} &=& 
\left( \matrix{ &  3y\pd{x} & 2x & -1 & 0 & \cr 
         &  2x\pd{y} & -3y\pd{y}+2 & -\pd{x} & 0 & \cr 
         &  -9y^2\pd{y}-3y & 0 & 2x & 1 & \cr 
         &  4x^2\pd{y} & 0 & -3y\pd{y}+4 & \pd{x} & \cr  }\right), \\
L^{-3} &=& 
\left( \matrix{ & -3y\pd{y}+2, & -2x, & -\pd{x}, & 1 & } \right).
\end{eqnarray*}
The shift vectors are given by
\begin{eqnarray*}
(s^1_1,s^1_2,s^1_3,s^1_4) &=& (0,1,0,-1), \\
(s^2_1,s^2_2,s^2_3,s^2_4) &=& (0,1,-1,0), \\
s^3_1 &=& 0.
\end{eqnarray*}
By computing the truncated complex, which is a complex of finite dimensional
vector spaces and linear maps, we obtain the result.

\end{example}


Programs written in the user language of {\tt kan/sm1} for algorithms 
in the present paper are available; please contact 
{\tt oaku@yokohama-cu.ac.jp}.

\section{Computation of cohomology groups with coefficients in a locally
constant sheaf of rank one}

A sheaf $\Vsc$ on $U$ is called {\it a locally constant sheaf} of rank
$m$
if for any $x \in U$, there exists an open set $W \ni x$ such that
the restriction $\Vsc_{|_W}$ is a constant sheaf $\C_W^m$.

Let $f_1, \ldots, f_d \in \Q[x]$ be (not necessarily irreducible) factors of
$f$ satisfying $f = f_1 \cdots f_d$.
Let $a_1, \ldots, a_d$ be complex numbers which lie in a computable
field.

The left $A_n$-module 
$$L(a) = \Q[x,1/f] f_1^{a_1} \cdots f_d^{a_d} $$
is defined as follows;
we define the action of $\pd{k}$ and $x_k$ by 
\begin{eqnarray*}
 \pd{k} \cdot ((g(x)/f^p) \cdot m ) &=& \left( \sum_{i=1}^d
     { {a_i {{\partial f_i} \over {\partial x_k}}} \over
       {f_i}} \right) (g(x)/f^p)  \cdot m 
     + {{\partial (g/f^p)} \over {\partial x_k}} \cdot m, \\
 x_k \cdot ((g(x)/f^p) \cdot m) &=& x_k g(x)/f^p \cdot m \\
\end{eqnarray*}
where $m = f_1^{a_1} \cdots f_d^{a_d}$ and
$g(x)$ is an arbitrary polynomial.
In fact, we can easily check that
$$ \pd{k} \cdot ( x_k (g/f^p) m) = ( \pd{k} x_k ) \cdot ((g/f^p) m) $$
and hence our definition of the action is well-defined.

The left $A_n$-module 
$$ P(a) = A_n f_1^{a_1} \cdots f_d^{a_d} $$ is 
the left $A_n$-submodule of $L(a)$ generated by $m$.

Put 
$$ \Vsc = \Hom_{A_n} (P(a), \hol_U^{an}), $$
where $\hol_U^{an}$ is the sheaf of {\it holomorphic} functions on
the complex manifold $U = \C^n \setminus V(f)$.
Here we endow $U$ with the classical topology instead of the Zariski topology.
When the left $A_n$-module $P(a)$ is expressed as
$ A_n/I(a)$,  we can regard ${\cal V}$ as a sheaf of holomorphic solutions
on $U$ of the system of linear partial differential equations $I(a)$;
we have, for a simply connected  open set $u \subset U_{cl}$,
$$ {\cal V}(u) \simeq \{ f \in \hol_U^{an}(u) \,|\,  \ell f = 0 \mbox{ for all }
  \ell \in I(a) \}
$$ 
where the isomorphism is given by
$$  {\cal V}(u) \ni \varphi \mapsto \varphi(1) \in \hol_U^{an}(u). $$
The $\C$-vector space ${\cal V}(u)$ is one dimensional and 
spanned by the function
$ f_1^{a_1} \cdots f_d^{a_d} $
since $I(a)$ contains
$$
\pd{k} - \sum_{i=1}^d{ {a_i {{\partial f_i} \over {\partial x_k}}} \over
       {f_i}}
\qquad (k = 1,\dots,n).
$$
Thus, ${\cal V}$ is the locally constant sheaf of rank one.

\begin{theorem}  \label{theorem:locallyConstantSheaf}
The cohomology group
$H^k(U,{\cal V})$ is computable.
\end{theorem}

This theorem is a generalization of Theorem \ref{theorem:constantSheaf}.
In fact, when $a_1 = \ldots = a_d = 0$,
the locally constant sheaf ${\cal V}$ is the constant sheaf $\C_U$.
In order to prove this theorem, we need to generalize Procedure
\ref{procedure:annihilator}
to compute a left ideal $I(a)$ of $A_n$ 
such that
$ L(a) = A_n/I(a)$
where $I(a)$ is, intuitively speaking,
the differential equations for
$(f^{-\nu})f_1^{a_1} \cdots f_d^{a_d}$ with an appropriate nonnegative
integer $\nu$.

We introduce the Weyl algebra
$$A_{d+n} = \Q\langle t_1, \ldots, t_d, x_1, \ldots, x_n,
  \pd{t_1}, \ldots, \pd{t_d}, \pd{1}, \ldots, \pd{n} \rangle.
$$
for our computation of $I(a)$.

\begin{procedure} \rm \label{procedure:severalAnnihilator}
(Computing $L(a)$.)
\\

\noindent
Input: $f, f_1, \ldots, f_d, a_1, \ldots, a_d$. \\
Output: a left ideal $I(a)$ of $A_n$ such that 
$L(a) = \Q[x,1/f]f_1^{a_1} \cdots f_d^{a_d} \simeq A_n/I(a)$.

\begin{enumerate}
\item (Computation of the annihilating ideal of $f_1^{s_1} \cdots f_d^{s_d}$ 
with indeterminates $s_1,\dots,s_d$)

Compute
\begin{eqnarray*}
& &\langle t_j-f_j(x) \ (j=1, \ldots, d), 
                {{\partial f_j} \over
                 {\partial x_i}} \pd{t_j} + \pd{i} \ 
                (i=1, \ldots, n, j=1, \ldots, d) \rangle \\
&\cap&
\Q[ t_1 \pd{t_1}, \cdots, t_d \pd{t_d}] \langle x, \pd{x} \rangle.
\end{eqnarray*}
Replacing each $t_i \pd{t_i}$ by the indeterminate $-s_i-1$ in generators
of the intersection, 
we obtain the set
$$G_0(-s_1-1, \ldots, -s_d-1) =
\{ Q_1(x,\pd{x},-s-1), \ldots, Q_k(x,\pd{x},-s-1) \}.
$$
(Call Procedure \ref{procedure:intersectionWithSubring} 
to compute the intersection
of a left ideal and the subring 
$\Q[t_1\pd{t_1}, \ldots, t_d \pd{t_d}]\langle x, \pd{x}\rangle$.)
The left ideal $I(s)$ of
$\Q[t_1\pd{t_1}, \ldots, t_d \pd{t_d}]\langle x, \pd{x}\rangle$
generated by $G_0(-s-1)$ gives the annihilating ideal for
$ f_1^{s_1} \cdots f_d^{s_d}$.

\item
Compute 
$$  \langle I(s), f_1(x), \ldots, f_d(x) \rangle \cap \Q[s_1, \ldots, s_d] $$
by an elimination order $ x, \pd{x} > s_1, \ldots, s_d  $.
Let $G_1(s)$ be a set of generators of the elimination ideal above.

\item
Choose a positive integer $\nu$ such that the set
$$(a_1-\nu, \ldots, a_d-\nu) - {\bf Z}_{>0}(1, \ldots, 1) $$
is not contained in the zero set
$$V(G_1(s)) = \{ v \in \C^d \,|\,
  g(v) = 0  \mbox{ for all }  g(s) \in G_1(s)
\}.
$$

\item
Output
$$I(a):=G_0(-s-1)_{|_{ s_1 \mapsto a_1-\nu, \ldots, s_d \mapsto a_d-\nu}}. $$
\end{enumerate}
\end{procedure}

In the above procedure,  $I(a)$ is the annihilating ideal of 
$f_1^{a_1-\nu}\cdots f_d^{a_d-\nu}$.  
The annihilating ideal of $f_1^{a_1}\cdots f_d^{a_d}$ can be computed
as the ideal quotient $I(a):(A_n f^\nu)$ through syzygy computation 
by means of Gr\"obner basis.     

Let us present an algorithm to compute the cohomology groups
 $H^k(U, {\cal V})$.

\begin{algorithm} \rm
\label{algorithm:cohomology-top2}
(Computing the cohomology groups $H^k(U,\Vsc)$.)
\\

\noindent
Input: $f, f_1, \ldots, f_d, a_1, \ldots, a_d$. \\
Output: the cohomology groups $H^k(U,\Vsc)$.

\begin{enumerate}
\item Call Procedure \ref{procedure:severalAnnihilator} with the input
$f, f_1, \ldots, f_d, -a_1, \ldots, -a_d$.
Get the output $I(-a)$.
\item Call Procedure \ref{procedure:integral} 
with the input $I = I(-a)$.
\end{enumerate}
\end{algorithm}

\section{Computation of $P(a)$ and its localization}

\setcounter{equation}{0}

Put $X = \C^n$ and let $Y$ be an algebraic set of $X$ defined by 
the polynomial $f \in \Q[x]$ with $x = (x_1,\dots,x_n)$. 
Let $\partial = (\partial_1,\dots,\partial_n)$ be the corresponding
differentiations.   
We denote by $\Osc_X$ and 
$\Dsc_X=\Osc_X \langle \pd{1}, \ldots, \pd{n} \rangle$ 
the sheaf of regular functions,  and
the sheaf of algebraic differential operators on $X$ respectively  
(see, e.g., \cite[p.15 and p.70]{HartshornBook} and \cite[p.15]{HottaTanisaki}).
We note that the set of the global sections
$\Gamma(X, \Dsc_X)$ coincides with $\C\otimes_{\Q}A_n$,
which is the Weyl algebra with coefficients in the complex numbers.
We will denote it also by $A_n$ if there is no risk of confusion.

In the sequel, we shall work in the category of algebraic $\Dsc_X$-modules
and prove isomorphisms for sheaves of $\Dsc_X$-modules.
Correctness of algorithms and procedures given in preceding sections follows 
by taking global section on $X$ in isomorphisms of propositions.

Put
$$
\Msc = \Psc(a) = 
 \Dsc_X f_1^{a_1}\cdots f_d^{a_d}.
$$
The left coherent $\Dsc_X$-module
$\Msc$ is a locally free $\Osc_X$-module of rank one
on $X\setminus Y$,
which is called an integrable connection and $\Msc_{|_{X \setminus Y}}$
has regular singularities along $Y$.    
As we have seen,
in order to compute the cohomology group, we may compute
$\int_X^{\cdot} \Msc[1/f]$.
Our purpose in this section is to give a proof of correctness of
Procedure \ref{procedure:severalAnnihilator}, which also gives
an algorithm to compute the localization 
$\Msc[1/f]:= \Osc_X[1/f]\otimes_{\Osc_X}\Msc$.
$\Msc[1/f]$ is 
a holonomic system on $X$ (Theorem 1.3 of Kashiwara \cite{K1}) 
and coincides with $\Msc$ on $X\setminus Y$.

We outline a method to compute $\Psc(a)[1/f]$ for given 
non-constant polynomials $f_1,\dots,f_d$ and $a = (a_1,\dots,a_d)$ 
with $f := f_1\cdots f_d$.  
Here, we assume that $a_i$ lies in a computable field.

Let $s = (s_1,\dots,s_d)$ be commutative indeterminates and put
$$ \Lsc(s) := \Osc_X[s,1/f]f_1^{s_1}\cdots f_d^{s_d},$$
which we regard as a free $\Osc_X[s,1/f]$-module.
Put $\Psc(s) := \Dsc_X[s]f_1^{s_1}\cdots f_d^{s_d}$.  
Then the set of global sections $\Gamma(X,\Lsc(s))$ of $\Lsc(s)$ coincides 
with $\C[x,s,1/f]f_1^{s_1}\cdots f_d^{s_d}$, and that of $\Psc(s)$ with
$ A_n[s]f_1^{s_1}\cdots f_d^{s_d}$. 

\begin{definition}{\rm
The (global) Bernstein-Sato ideal $B(f_1,\dots,f_d)$ of $\Q[s]$ is defined by
$$
B(f_1,\dots,f_d) := \{ b(s) \in \Q[s] \mid 
b(s)f_1^{s_1}\cdots f_d^{s_d} \in 
A_n[s]f_1^{s_1+1}\cdots f_d^{s_d+1}\}.
$$
}
\end{definition}

The step 2 of Procedure \ref{procedure:severalAnnihilator} gives
an algorithm to compute the Bernstein-Sato ideal.

\begin{proposition}[\cite{S}]
\label{proposition:sabbahBfunction}
There exist a finite number of linear forms $L_1(s),\dots,L_\kappa(s)$ 
in $s$ with nonnegative integer coefficients, and nonzero 
univariate polynomials $b_1,\dots,b_\kappa$, such that
$$ 
b(s) := b_1(L_1(s))\cdots b_\kappa(L_\kappa(s)) \in B(f_1,\dots,f_d).
$$
In particular, for any $a = (a_1,\dots,a_d) \in \C^d$, 
the intersection of $\{(a_1-\nu,\dots,a_d-\nu) \mid \nu\in\N\}$ with
$$
\V(B(f_1,\dots,f_d)) := \{s = (s_1,\dots,s_d) \in \C^d \mid
b(s) = 0 \mbox{ for any $b \in B(f_1,\dots,f_d)$}\}
$$
is a finite set.  
\end{proposition}

The following proposition tells us that
if $a$ is generic, then the localization $\Lsc(a)$ of $\Psc(a)$
agrees with $\Psc(a)$.

\begin{proposition}  \label{proposition:generic}
Assume that $a = (a_1,\dots,a_d) \in \C^d$ satisfy that 
$(a_1-\nu,\dots,a_d-\nu)$ is not contained in $\V(B(f_1,\dots,f_d))$ 
for any $\nu = 1,2,3,\dots$.
Then $\Psc(a) = \Lsc(a)$ holds.
In particular, the $\Osc_X$-homomorphism 
$f : \Psc(a) \longrightarrow \Psc(a)$ is an isomorphism.	
\end{proposition}

\begin{proof}
In the notation of Proposition \ref{proposition:sabbahBfunction}, there exist
$b(s)\in B(f_1,\dots,f_d)$ and $p(s)\in A_n[s]$ such that
$$
p(s)f_1^{s_1+1}\cdots f_d^{s_d+1} = b(s)f_1^{s_1}\cdots f_d^{s_d}.
$$ 
and $b(a_1-1,\dots,a_d-1)\neq 0$. Then we have
$$
f_1^{a_1-1}\cdots f_d^{a_d-1} = b(a_1-1,\cdots,a_d-1)^{-1}
p(a_1-1,\dots,a_d-1)f_1^{a_1}\cdots f_d^{a_d}.
$$
Proceeding in the same way by using the assumption, we know that
$f_1^{a_1-\nu}\cdots f_d^{a_d-\nu}$ is contained in $\Psc(a)$ for 
$\nu = 1,2,3,\dots$.  This implies $\Psc(a) = \Lsc(a)$.
\end{proof}

Next, we shall see that the localization $\Lsc(a)$ agrees with 
$\Psc(a_1-\nu_0, \ldots, a_d-\nu_0)$ for an integer $\nu_0$
determined by the zero set of the Berndstein-Sato ideal.
In order to prove this fact, we need a lemma.

\begin{lemma}
$\Osc_X[1/f]$ is a flat $\Osc_X$-module.
\end{lemma}

\begin{proof}
This should be well-known (e.g. this is a special case of 
Lemma 1.1 of \cite{KK}).  Here we give a direct proof.
Let $\iota : \Ksc \longrightarrow \Nsc$ be an arbitrary injective 
$\Osc_X$-homomorphism.  
Then for a section $u$ of $\Ksc$, we have $1\otimes u = 0$ in 
$\Ksc[1/f] = \Osc_X[1/f]\otimes_{\Osc_X}\Ksc$ if and only if 
$f^\nu u = 0$ for some $\nu \in \N$ (cf.\ Lemma 7.2 of \cite{OakuAdvance}). 
Hence we have $1\otimes u = 0$ if and only if $1\otimes \iota(u) = 0$.  
This completes the proof.
\end{proof}

\begin{proposition}  \label{proposition:3.5}
Fix an arbitrary $a = (a_1,\dots,a_d) \in \C^d$.
Let $\nu_0$ be a positive integer such that 
$(a_1-\nu,\dots,a_d-\nu)$ is not contained in $\V(B(f_1,\dots,f_d))$ 
for any integer $\nu > \nu_0 $.  
Then we have
$$
\Psc(a)[1/f] = \Lsc(a) = \Psc(a_1-\nu_0,\dots,a_d-\nu_0). 
$$
\end{proposition}

\begin{proof}
Consider the short exact sequence
\begin{equation}
\label{equation:3.1}
0 \longrightarrow \Psc(a) \stackrel{\iota}{\longrightarrow} \Lsc(a) 
\longrightarrow
\Lsc(a)/\Psc(a) \longrightarrow 0,  
\end{equation}
where $\iota$ is the inclusion.  
First note that $(\Lsc(a)/\Psc(a))[1/f] = 0$.
In fact, any section $v$ of $\Lsc(a)$ is written in the form
$v = g f_1^{a_1-\nu}\dots f_d^{a_d-\nu}$ with $g \in \Osc_X$ and $\nu \in\N$.
Hence we have $f^\nu v \in \Psc(a)$.  
This implies $(\Lsc(a)/\Psc(a))[1/f] = 0$.

Since $\Osc_X[1/f]$ is a flat $\Osc_X$-module, 
we have from (\ref{equation:3.1}) an exact sequence 
$$ 
0 \longrightarrow \Psc(a)[1/f] 
\stackrel{1\otimes\iota} \longrightarrow \Lsc(a)[1/f] 
\longrightarrow 0.  
$$
Since $\Lsc(a)[1/f] = \Lsc(a)$, we have proved
the first equality of the proposition.
The second one follows from Proposition \ref{proposition:generic} since
$\Lsc(a) = \Lsc(a_1-\nu_0,\dots,a_d-\nu_0)$ 
($f$ is invertible in $\Lsc(a)$).
\end{proof}

\begin{proposition} \label{proposition:3.6}
Under the same assumption as Proposition \ref{proposition:generic}, 
the 
$\Dsc_X$-homomorphism (specialization $s=a$)
$$
\rho :\,\, \Psc(s)/((s_1-a_1)\Psc(s) + \dots + (s_d-a_d)\Psc(s)) 
\longrightarrow \Psc(a)
$$ 
is an isomorphism.  
\end{proposition}

\begin{proof}
Assume that a section $u := p(s)f_1^{s_1}\cdots f_d^{s_d}$ of $\Psc(s)$
satisfies $\rho(\overline u) = 0$, where $\overline u$ denotes the modulo
class of $u$. 
Then there exist $g_1(s),\dots,g_d(s) \in \Osc_X[s]$ and $\nu \in \N$ 
such that
$$
u = \sum_{j=1}^d (s_j-a_j)g_j(s) f_1^{s_1-\nu}\cdots f_d^{s_d-\nu}.
$$
By the same argument as the proof of Proposition \ref{proposition:generic},
we can find $\tilde b(s) \in \Q[s]$ and $Q(s) \in \Dsc_X[s]$ such that
$$ 
\tilde b(s) f_1^{s_1-\nu}\cdots f_d^{s_d-\nu} 
= Q(s)f_1^{s_1}\cdots f_d^{s_d}
$$
and $\tilde b(a) \neq 0$.

There exist $c_1(s),\dots,c_d(s) \in \C[s]$ which satisfy
$$ \tilde b(a) - \tilde b(s) = \sum_{j=1}^d (s_j-a_j)c_j(s).$$
Hence we get
\begin{eqnarray*}
\tilde b(a) u &=& 
\left(\tilde b(s)p(s) + \sum_{j=1}^d (s_j-a_j)c_j(s)p(s)\right)
f_1^{s_1}\cdots f_d^{s_d}\\
&=&
\sum_{j=1}^d (s_j-a_j) \left( \tilde b(s) g_j(s)
f_1^{s_1-\nu}\cdots f_d^{s_d-\nu} 
+ c_j(s)p(s)f_1^{s_1}\cdots f_d^{s_d} \right) \\
&=&
\sum_{j=1}^d (s_j-a_j)\left(g_j(s)Q(s) + c_j(s)p(s)\right) 
f_1^{s_1}\cdots f_d^{s_d}.
\end{eqnarray*}
Since $\tilde b(a) \neq 0$ by the assumption, we conclude that
$u \in (s_1-a_1)\Psc(s) + \cdots + (s_d-a_d)\Psc(s)$.  
Hence $\rho$ is injective. The surjectivity is obvious.
\end{proof}

Let us consider the problem of finding the annihilating ideal
of $f_1^{s_1} \cdots f_d^{s_d}$.

Let $A_d$ be the Weyl algebra on the variables $t = (t_1,\dots,t_d)$.  
We denote by $A_d\Dsc_X := A_d \otimes_\C\Dsc_X$ 
the sheaf on $X$ of the
differential operators in variables $(t,x)$ which are polynomials in $t$.
We follow an argument of Malgrange \cite{Malgrange} for the case of $d=1$.

We can endow $\Lsc(s)$ with a structure of left $A_d\Dsc_X$-module by
\begin{eqnarray}
& & t_j(g(x,s)f_1^{s_1}\cdots f_d^{s_d})  \nonumber \\
&=& 
g(x,s_1,\dots,s_j+1,\dots,s_d)f_1^{s_1}\cdots f_j^{s_j+1}\cdots f_d^{s_d},
\label{equation:module1}
\\
& & \partial_{t_j}(g(x,s)f_1^{s_1}\cdots f_d^{s_d}) \nonumber \\
&=& -s_j
g(x,s_1,\dots,s_j-1,\dots,s_d)f_1^{s_1}\cdots f_j^{s_j-1}\cdots f_d^{s_d}
\label{equation:module2}
\end{eqnarray}
for $g(x,s) \in \Osc_X[s,1/f]$ and $j=1,\dots,d$.  

\begin{lemma}
Let $\Nsc$ be a sheaf of left ideals of $A_d\Dsc_X$ generated by
\begin{eqnarray}
t_j - f_j(x) &\quad& (j=1,\dots,d),
\label{equation:annihilator1}\\
\partial_{x_i} + \sum_{j=1}^d\frac{\partial f_j}{\partial x_i}\partial_{t_j}
&\quad& (i=1,\dots,n).
\label{equation:annihilator2}
\end{eqnarray}
Then each stalk of $\Nsc$ is a maximal left ideal.
\end{lemma}

\begin{proof}
By a coordinate transformation 
$$
t'_j = t_j - f_j(x) \quad (j = 1,\dots,d),\qquad
x' = x,
$$
we can reduce to the case where $f_1 = \dots = f_d = 0$.  
In that case, the statement is obvious.
\end{proof}

\begin{proposition}
We have
$$
\Nsc = \{p \in A_d\Dsc_X \mid 
p f_1^{s_1}\cdots f_d^{s_d} = 0\}.
$$
\end{proposition}

\begin{proof}
It is easy to verify the inclusion $\subset$ by using 
(\ref{equation:module1}) and (\ref{equation:module2}).  
Since $\Nsc$ is maximal, we obtain the equality.
\end{proof}

We put
$$
\Isc(s) := \{ p(s) \in \Dsc_X[s] \mid
p(s)f_1^{s_1}\cdots f_d^{s_d} = 0\}.
$$

\begin{proposition} \label{proposition:severalAnnihilator}
For a Zariski open set $u$ of $X$, we have
\begin{eqnarray*}
& & \Gamma(u,\Isc(s)) \\
&=& \{p(-s_1-1,\dots,-s_d-1) \mid 
p(t_1\partial_{t_1},\dots,t_d\partial_{t_d}) \in 
\Gamma(u, \Nsc\cap\Dsc_X[t_1\partial_{t_1},\dots,t_d\partial_{t_d}]) \}.
\end{eqnarray*}
\end{proposition}

\begin{proof}
By (\ref{equation:annihilator1}) and 
(\ref{equation:annihilator2}), we get the relations
$$
s_j = -\partial_{t_j}t_j = -t_j\partial_{t_j}-1
\quad (j=1,\dots,d).
$$
Hence $\Dsc_X[s]$ is isomorphic to the subring
$\Dsc_X[t_1\partial_{t_1},\dots,t_d\partial_{t_d}]$ of $A_d\Dsc_X$.
This implies the conclusion.
\end{proof}

\begin{proposition}
\label{proposition:3.10}
Procedure
\ref{procedure:severalAnnihilator} is correct.
\end{proposition}

\begin{proof}
The correctness of step 1 follows from Proposition 
\ref{proposition:severalAnnihilator}.

To verify the correctness of step 2 of Procedcure
\ref{procedure:severalAnnihilator},
one has only to note that for $b(s) \in \Q[s]$, we have
$b(s)\in B(f_1,\dots,f_d)$ if and only if 
$b(s)-f$ belongs to $\Gamma(X,\Isc(s))$.

The correctness of steps 3 and 4 can be shown 
by taking global sections in sheaf isomorphisms given
in Propositions \ref{proposition:3.5} and \ref{proposition:3.6}.
\end{proof}

As to our experiments, it is more efficient that one eliminates 
$\partial_x$ first, and then elminates $x$ in step 2 of 
Procedure \ref{procedure:severalAnnihilator}.
However, even with this, the complexity of Procedure 
\ref{procedure:severalAnnihilator} is huge.

\section{Computation of the intersection of a left ideal and a subring}
\setcounter{equation}{0}
\label{section:intersection}

In this section, we give a procedure to compute
the intersection of the left ideal
$$
\langle t_j-f_j(x) \ (j=1, \ldots, d), 
                {{\partial f_j} \over
                 {\partial x_i}} \pd{t_j} + \pd{i} \ 
                (i=1, \ldots, n, j=1, \ldots, d)
\rangle $$ 
in $A_{d+n}$ and the subring
$ \Q[ t_1 \pd{t_1}, \cdots, t_d \pd{t_d}] \langle x, \pd{x} \rangle $
of $A_{d+n}$.
The intersection gives the annihilating ideal for
$f_1^{s_1} \cdots f_d^{s_d}$ with the replacement
$ t_i \pd{t_i} \mapsto -s_i-1$.

\begin{procedure}
\label{procedure:intersectionWithSubring}
{\rm

\noindent
Input: polynomials $f_1,\dots,f_d$ in $x = (x_1,\dots,x_n)$. \\
Output:  a set of generators of the annihilating ideal
$\Isc(s)$
of $f_1^{s_1} \cdots f_d^{s_d}$.

\begin{enumerate}
\item
Introducing indeterminates $t = (t_1,\dots,t_d)$, $u = (u_1,\dots,u_d)$, 
$v = (v_1,\dots,v_d)$, let $I$ be the left ideal of 
$A_{n+d}[u,v] = \Q[u,v]\langle x,t, \partial_x,\partial_t\rangle$
generated by
\begin{eqnarray}
t_j-u_jf_j,\,\, &\quad& (j=1,\dots,d),
\label{equation:uv1}\\
\partial_{x_i} + \sum_{j=1}^d\frac{\partial f_j}{\partial x_i}
u_j\partial_{t_j} &\quad& (i=1,\dots,n),
\label{equation:uv2}\\
1-u_jv_j,\,\, &\quad& (j=1,\dots,d).
\label{equation:uv3}
\end{eqnarray}
\item
Take any term order on $A_{n+d}[u,v]$ for eliminating 
$u,v$.  Let $G$ be a Gr\"obner basis of $I$ with respect to this term order.
Put $G_0 = \{P_1,\dots,P_k\} := G \cap A_{n+d}$.  
\item
For each $i = 1,\dots,k$, there exist $Q_i \in \Dsc_X[s]$ and 
$\nu_{i1},\dots,\nu_{id} \in \Z$ such that
$$
S_{1,\nu_{i1}}\cdots S_{d,\nu_{id}}P_i = 
Q_i(x,\partial_x, t_1\partial_{t_1},\dots,t_d\partial_{t_d}) 
$$ 
holds, where $S_{j,\nu} := \partial_{t_j}^\nu$ if $\nu \geq 0$, 
and $S_{j,\nu} := t_j^{-\nu}$ otherwise.  Set
$$
G_0(s) := \{Q_1(x,\partial_x,s),\dots,Q_k(x,\partial_x,s)\}.
$$
\end{enumerate}
Output: $G_0(-s_1-1,\dots,-s_d-1)$ is a set of generators of $\Isc(s)$.
}
\end{procedure}

\begin{proposition}
\label{proposition:4.2}
Procedure \ref{procedure:intersectionWithSubring} is correct.
\end{proposition}

\begin{proof}
First, we must show that each element of $G_0$ can be written in the form 
as in the step 3 of Procedure \ref{procedure:intersectionWithSubring}. 
Fix any $j$ with $1 \leq j \leq d$.  Then
the generators of $I$ given in the step 1 are
homogeneous with respect to the weight table ${\cal W}_j$ below: 
\begin{center} ${\cal W}_j$: 
\begin{tabular}{|l||c|c|c|c|c|c|} \hline
variables & $x_i,\partial_{x_i}\,\,(1\leq i \leq n)$ & $t_j$ & $\partial_{t_j}$
& $u_j$ & $v_j$ & $t_k,\partial_{t_k},u_k,v_k\,\,(k \neq j)$ \\ \hline
weight  & $0$ & $-1$ & $1$ & $-1$ & $1$ & $0$ \\ \hline
\end{tabular}
\end{center}
Moreover, the product of two operators preserves the homogeneity with respect 
to ${\cal W}_j$.
Hence each element of $G_0$ is homogeneous with respect to ${\cal W}_j$ and
free of $u$ and $v$.  This enables us to write $P_i$ in the form as in 
the step 3.

Now let us show that each $Q_i(x,\partial_x,-s-1)$ belongs to $\Isc(s)$
with the notation $-s-1 = (-s_1-1,\dots,-s_d-1)$.  
By the definition, $P_i$ is contained in the ideal generated by
(\ref{equation:uv1})--(\ref{equation:uv3}).  
Substituting 1 for every $u_i$ and $v_i$, we know that
$P_i$ belongs to $\Nsc$ since it does not depend on $u,v$.
Hence $Q_i(-s-1)$  belongs to $\Isc(s)$ in view of 
Proposition \ref{proposition:severalAnnihilator}.

Conversely, let $p(-s-1)$ be an arbitrary section of $\Isc(s)$.
Multiplying by a polynomial, we may assume that 
$p(t_1\partial_{t_1},\dots,t_d\partial_{t_d})$ belongs to the left 
ideal of $A_{n+d}$ generated by (\ref{equation:annihilator1}) and 
(\ref{equation:annihilator2}) making use of 
Proposition \ref{proposition:severalAnnihilator}
again.
That is, there exist $R_j,S_i\in A_{n+d}$ such that
\begin{equation} \label{equation:afo1}
p(t_1\partial_{t_1},\dots,t_d\partial_{t_d})
= \sum_{j=1}^d R_j\cdot(t_j - f_j)
+ \sum_{i=1}^n S_i\cdot
\left(\partial_{x_i} + \sum_{j=1}^d\frac{\partial f_j}{\partial x_i}\right)
\end{equation}
We can homogenize the both sides of 
(\ref{equation:afo1}) by adding $u$ 
with respect to the weight table ${\cal W}_j$.  
By performing this procedure for every $j=1,\dots,d$, we obtain 
a homogenization of 
(\ref{equation:afo1}) with respect to all ${\cal W}_1,\dots,{\cal W}_d$.  
The left hand side of this homogenization is in the form 
$u_1^{\mu_1}\cdots u_d^{\mu_d}p$ with nonnegative integers 
$\mu_1,\dots,\mu_d$ since $p$ itself is homogenous.  
Thus $u_1^{\mu_1}\cdots u_d^{\mu_d}p$ is 
contained in the ideal of $A_n[u]$ generated by 
(\ref{equation:uv1}) and (\ref{equation:uv2}).  
This implies that
$$
p = (1-u_1^{\mu_1}\cdots u_d^{\mu_d}v_1^{\mu_1}\cdots v_d^{\mu_d})p
    + u_1^{\mu_1}\cdots u_d^{\mu_d}v_1^{\mu_1}\cdots v_d^{\mu_d}p
$$
belongs to $I$.  Since $G$ is a Gr\"obner basis of $I$ with respect to a term order
for eliminating $u,v$, there exist $U_1,\dots,U_k \in A_{n+d}$ such that
$$
p(t_1\partial_{t_1},\dots,t_d\partial_{t_d}) 
= \sum_{i=1}^k U_iP_i.
$$
Since $p$ and $P_i$ are homogeneous with respect to each ${\cal W}_j$, 
we may assume that so is $U_i$.  
Moreover, since the weight of $p$ is zero with respect to each ${\cal W}_j$, 
all $U_i$ are written in the form
$$
U_i = U'_i(t_1\partial_{t_1},\dots,t_d\partial_{t_d})
      S_{1,\nu_{i1}}\cdots S_{d,\nu_{id}}
$$
with some $U'_i \in A_n[t_1\partial_{t_1},\dots,t_d\partial_{t_d}]$.
Hence $p(s)$ belongs to the left ideal of $A_n[s]$ generated by
$G_0(s)$.  
This completes the proof.
\end{proof}

\begin{example}  \rm \label{example:example4.1}
Consider
$f = x^{s_1} y^{s_2} (1-x-y)^{s_3}$.
$\Isc(s)$ is generated by
\begin{eqnarray*}
 & &y s_1+y s_2+y s_3-s_2-y x \pd{x}-y^2 \pd{y}+y \pd{y} ,  \\
 & &x s_1+x s_2+x s_3-s_1-x^2 \pd{x}-y x \pd{y}+x \pd{x} , \\ 
 & & x s_2+y s_2+y s_3-s_2-y x \pd{y}-y^2 \pd{y}+y \pd{y}.
\end{eqnarray*}
Note that this ideal is strictly larger than the ideal generated by
trivial annihilators
$$ x(1-x-y) \pd{x} - x (1-x-y) (\partial f/\partial x)/f,\ 
   y(1-x-y) \pd{y} - y (1-x-y) (\partial f/\partial y)/f.
$$
\end{example}

\section{Twisted de Rham cohomology group}
\label{deRhamCohomologyAndIntegrationOfDmodules}
\setcounter{equation}{0}

In this section, we shall explain that computation of $\Dsc$-module
theoretic integrals  of
$L(a)$ gives the cohomology groups 
$H^k(U,{\cal V})$,
which is nothing but what Grothendieck-Deligne comparison 
theorem says;
the contents of this section should be well-known to specialists.
However, they are not explicitly explained in literatures.

First let us recall the integration functor for $\Dsc$-modules.
In general, let $\Msc$ be a left $\Dsc_X$-module 
(or, more generally, a complex of $\Dsc_X$-modules) defined on $X$.
Then integration of $\Msc$ over $X$ is defined by
$$
\int_X \Msc := R\Gamma(X, \Omega_X\otimes_{\Dsc_X}^{L}\Msc)
$$
as an object of the derived category of $\C$-vector spaces, 
where $R$ and $L$ denotes the right and the left derived functors in 
the derived categories, $\Gamma$ is the global section functor,
 and $\Omega_X$ is the sheaf of algebraic 
$n$-forms on $X$, which has a natural structure of the right $\Dsc_X$-module
and is isomorphic to
$\Dsc_X/ (\pd{1} \Dsc_X + \cdots + \pd{n} \Dsc_X)$  
since $X$ is the affine space. 
For $i\in \Z$, the $i$-th cohomology of $\int_X\Msc$ is denoted by 
$\int_X^i\Msc$, which is a $\C$-vector space.  
$R^i\Gamma(X,\Nsc)$ is often denoted by $H^i(X,\Nsc)$.
See, e.g.,  \cite{EisenbudBook} and \cite{HartshornBook}
for an introduction to the mechanism of derived functors.

Now put
$$
h_i = \sum_{j=1}^d a_j\frac{f}{f_j}\frac{\partial f_j}{\partial x_i}
\quad (i = 1,\dots,n).
$$
Let 
$\Msc$ be the left $\Dsc_X$-module $\Msc := \Dsc_X/\Isc$, where
$\Isc$ is the left ideal generated by 
$f\partial_i-h_i$ ($i=1,\dots,n$) with polynomials $h_i$.
Here, we note that $h_i$ satisfy the
integrability condition
\begin{equation}
\frac{\partial}{\partial x_j}(h_i/f) =
\frac{\partial}{\partial x_i}(h_j/f)
\qquad 1 \leq i,j \leq n,
\end{equation}
and the function $f_1^{a_1} \cdots f_d^{a_d}$ is annihilated 
by the operators $f \pd{i} - h_i$.
$\Msc$ has regular singularities along 
(the non-singular locus of) $Y$ and also along 
the hyperplane at infinity of the projective space ${\bf P}^n$. 
$\Msc$  and $\Psc(a)$ are isomorphic as $\Dsc_X$-modules 
on $X\setminus Y$.  In fact, both are simple holonomic systems 
and there exists a natural $\Dsc_X$-homomorphism of $\Msc$ to $\Psc(a)$ 
which sends the modulo class of $1 \in \Dsc_X$ to 
$f_1^{a_1}\cdots f_d^{a_d}$.  
However, these two modules are not isomorphic on $X$ in general.

The de Rham complex ${\rm DR}(\Msc[1/f])$ of the localization 
$\Msc[1/f] := \Osc_X[1/f]\otimes_{\Osc_X}\Msc$ is defined by
\begin{equation}
\label{equation:rationalDeRham}
0 \longrightarrow \Omega_X^{0}\otimes_{\Osc_X}\Msc[1/f]
\stackrel{d}{\longrightarrow} \Omega_X^{1}\otimes_{\Osc_X}\Msc[1/f] 
\stackrel{d}{\longrightarrow}
\cdots
\stackrel{d}{\longrightarrow} \Omega_X^{n}\otimes_{\Osc_X}\Msc[1/f]
\longrightarrow 0.
\end{equation}
where $d$ is given by
$$
d(dx_{k_1}\wedge \cdots \wedge dx_{k_i}\otimes u)
= \sum_{j=1}^n dx_j\wedge dx_{k_1}\wedge \cdots \wedge dx_{k_i}
\otimes (\partial_j u)
$$
for $u \in \Msc[1/f]$. 
As $\Dsc_X[1/f]$-module (not as $\Dsc_X$-module!), 
there is an isomorphism 
$$ 
\Msc[1/f] \simeq \Dsc_X[1/f]/
(\Dsc_X[1/f](\partial_1 - h_1f^{-1}) + \dots 
+ \Dsc_X[1/f](\partial_n -h_nf^{-1})).
$$
Let $P$ be a section of $\Msc[1/f]$.  Then there exist $Q_i\in\Dsc_X[1/f]$
and $r \in \Osc_X[1/f]$ such that
$$
P = \sum_{i=1}^n Q_i(\partial_i - h_if^{-1}) + r.
$$
Such $r$ is determined uniquely. 
Then we define $\varphi(P) = r$.
Hence
$$ 
\varphi : \Msc[1/f] \longrightarrow \Osc_X[1/f]
$$
defines an isomorphism as $\Osc_X[1/f]$-module.  
By transforming the complex (\ref{equation:rationalDeRham}) by means of this $\varphi$, we get
the following complex that is isomorphic to (\ref{equation:rationalDeRham}):
\begin{equation}
\label{equation:5.3}
0 \longrightarrow \Omega_X^{0}[1/f]
\stackrel{\nabla}{\longrightarrow} \Omega_X^{1}[1/f] 
\stackrel{\nabla}{\longrightarrow}
\cdots
\stackrel{\nabla}{\longrightarrow} \Omega_X^{n}[1/f]
\longrightarrow 0,
\end{equation}
where $\nabla$, which is called the integrable connection,
is defined by 
$$
\nabla(u dx_{i_1}\wedge \cdots \wedge dx_{i_k})
= \sum_{j=1}^n \left(\frac{\partial u}{\partial x_j}+\frac{h_j}{f}u\right) 
dx_j\wedge dx_{i_1}\wedge \cdots \wedge dx_{i_k}
$$
for $u \in \Osc_X[1/f]$; in fact, we have 
$$
\partial_j\cdot u = u\partial_j + \frac{\partial u}{\partial x_j} 
\equiv u\frac{h_j}{f} + \frac{\partial u}{\partial x_j} 
$$
modulo $\Dsc_X[1/f](\partial_1-h_1f^{-1}) +\dots+ 
\Dsc_X[1/f](\partial_n-h_nf^{-1})$.
Thus the integral $\int_X \Msc[1/f] = R\Gamma(X, (\ref{equation:5.3}))$ 
is isomorphic to the complex
\begin{equation} \label{equation:rationalNablaComplex}
0 \longrightarrow \Gamma(X;\Omega_X^{0}[1/f])
\stackrel{\nabla}{\longrightarrow} \Gamma(X;\Omega_X^{1}[1/f]) 
\stackrel{\nabla}{\longrightarrow}
\cdots
\stackrel{\nabla}{\longrightarrow} \Gamma(X;\Omega_X^{n}[1/f])
\longrightarrow 0
\end{equation}
since $\Omega_X^{k}[1/f]$ is a quasi-coherent $\Osc_X$-module,
$X$ is affine and hence
$ H^k(X,\Omega_X^p[1/f]) = 0$ for $k > 0$
(see, e.g., \cite[p. 205, Proposition 1.2A, p.215, Theorem 3.7]{HartshornBook}
and \cite{Serre}).
The cohomology of this complex is nothing but the algebraic 
twisted de Rham cohomology with respect to the local system on 
$X\setminus Y$ defined by
the equation $\nabla u = 0$ for $u \in \Osc_X$.
When $\Msc = \Psc(a)$ on $X \setminus Y$, 
(\ref{equation:rationalNablaComplex}) 
gives the algebraic twisted de Rham 
cohomology groups associated with the local system defined by $\Psc(-a)$, 
i.e. the cohomology groups of $X\setminus Y$ 
with coefficients in the locally constant sheaf 
$$
{\cal V} := \Hom_{\Dsc_X}(\Psc(-a),\Osc_X^{\rm an})
= \Hom_{\Dsc_X}(\Osc_X, \Dsc_X^{\rm an}\otimes_{\Dsc_X}\Psc(a)),
$$
where $\Dsc^{\rm an}_X$ denotes the sheaf of holomorphic differential 
operators.  
In fact, by applying the functor $\Osc^{\rm an}_X\otimes_{\Osc_X}$ to
the complex (5.3), we obtain a complex of sheaves on $X\setminus Y$ 
whose $k$-th cohomology group is
$$
\{ u \in \Osc^{\rm an}_X \mid \nabla u = 0\} 
= \Hom_{\Dsc_X}(\Psc(-a),\Osc_X^{\rm an})
$$
if $k=0$ and zero otherwise.

The algebraic twisted de Rham cohomology coincides with the analytic
one by virtue of the comparison theorem of Deligne 
\cite[p.98 Theorem 6.2, p.99 Corollary 6.3]{D}.
Let us summarize what we have explained.

\begin{theorem} {\rm (Comparison theorem, \cite{D})}
\label{theorem:comparison}
$$H^k(U,\Vsc) 
\simeq H^{k}((\Gamma(X;\Omega_X^{\bullet}[1/f]),\nabla))
\simeq  H^{k-n}(X, DR(\Msc[1/f])).
$$
\end{theorem}


As we will see in Proposition \ref{proposition:6.1},
we, moreover, have
$$
H^{k-n}(X, DR(\Msc[1/f])) 
\simeq 
 H^{k-n}(A_n/(\pd{1} A_n + \cdots + \pd{n} A_n) \otimes_{A_n}^{ L}
 A_n/(I(a))).
$$

\begin{example} \rm (Beta function) 

Putting $X = \C$, we consider
$ \Psc(a) = \Dsc_X x^{a_1} (1-x)^{a_2} $
for generic complex numbers $a_1$ and $a_2$.
We have
$\Psc(a) = \Lsc(a) \simeq \Dsc_X/\langle p \rangle$
with $p = (x^2-x)\pd{x} - (a_1+a_2) x + a_1$.
The Bernstein-Sato ideal for $x$ and $1-x$ is generated by $(s_1+1)(s_2+1)$.  
The $b$-function of the Fourier transform $\Dsc_X/\langle \hat p \rangle$
with 
$$
 \hat p = x\pd{x}^2 + (x+a_1+a_2+2)\pd{x} + a_1 + 1
$$
is $s(s+a_1+a_2+1)$.  
Hence by applying Procedure \ref{procedure:integral} with $k_1 = 0$, 
we have
$$  0 \longrightarrow F_{-1}/( F_{-1} + x A_1 )
    \stackrel{\cdot {\hat p}}{\longrightarrow}
      F_0/( F_{-1} + x A_1 )
    \longrightarrow 0
$$ and
we get
$$ H^0(U,\Vsc) = 0, \quad H^1(U,\Vsc) = \C $$
where  $U = \{ 0, 1 \}$
and
$$\Vsc(w) =  \{ u \in \hol^{an}(w) \,|\, 
               du/dx = (-a_1/x+a_2/(1-x)) u \}$$
for a simply connected open set $w$.
Note that
$$ H^{1}(U,\Vsc) \simeq{{\C\left[x,{1 \over {x(1-x)}}\right]dx} \over
                        {\nabla \C\left[x,{1 \over {x(1-x)}}\right]}}
                 \simeq \C \cdot \left({1 \over x}-{1 \over {1-x}} \right) dx
$$
where $\nabla = d + (a_1/x - a_2/(1-x))dx \wedge$.
The beta function should be regarded as
$$ \int_{0}^1 x^{a_1}(1-x)^{a_2} \varphi $$
where
$\varphi={{dx} \over {x(1-x)}} \in H^{1}(U,\Vsc)$.
\end{example}

\begin{example} \rm
For generic complex numbers $a_1, \ldots, a_m$, 
we consider
$\Psc(a) = \Dsc_X \prod_{i=1}^m (x-c_i)^{a_i}$
where
$ c_1, \ldots, c_m$ are distinct points in $\C$.
By applying our algorithm, we can see that
$ H^1(U,\Vsc) = \C^{m-1}$ and
$H^0(U,\Vsc) = 0$
where $U = \C \setminus \{ c_1, \ldots, c_m \}$
and
$ \Vsc = \Hom_{\Dsc_U}(\Psc(-a),\hol^{an}_U)$.
See \cite{Aomoto-Kita} for details on these cohomology groups 
and hypergeometric functions.
\end{example}

\begin{example} \rm 
(Counting the number of bounded chambers by $\Dsc$-module algorithms)

We consider a collection of hyperplanes
$$ L_i(x) = \sum_{j=1}^n c_{ij} x_j + c_{i0} = 0, \quad (i=1, \ldots, m)$$
in $\R^n$ and put $f = \prod_{i=1}^m L_i(x)$.
For complex numbers $a_1, \ldots, a_m$, we consider
$ \Psc(a) = \Dsc_X \prod_{i=1}^m L_i(x)^{a_i} $.
The number of bounded chambers in 
$U = \R^n \setminus \cup_{i=1}^m \{ x \,|\, L_i(x) = 0 \}$
is equal to the Euler number of $H^k(U,\Vsc)$ 
(see \cite[p.47 Theorem 2.13.1]{Aomoto-Kita} and \cite{Orlik-Terao}).
Although there are several algorithms in computational geometry
to count the number,
this number can also be counted by our purely algebraic algorithm.
Table \ref{table:1} is an example of computation of Euler numbers 
by our algorithm and implementation.

\begin{table}[htb]
\caption{${\rm dim}_\C H^i(X\setminus Y, \C)$}
\label{table:1}
\begin{center}
\begin{tabular}{|c|c|c|c|c|}\hline
$f$ & $i=2$ & $i=1$ & $i=0$ & Euler ch. \\ \hline\hline
$xy$               & 1 & 2 & 1 & 0 \\ \hline
$xy(x+y+1)$        & 3 & 3 & 1 & 1 \\ \hline
\parbox{38mm}{$xy$\\$\cdot(x+y+1)$\\$\cdot(x-y-2)$} 
                   & 6 & 4 & 1 & 3 \\ \hline
\end{tabular}
\end{center}
where
$X = \C^2 = \{(x,y)\}$, $Y=\{f=0\}$.
\end{table}

Of course, our method is far from efficient.
However, it is rather surprising that purely algebraic 
computations in the ring of differential operators can 
evaluate the number of bounded chambers in a given hyperplane arrangement.
\end{example}

\section{Computation of integration}
\setcounter{equation}{0}
\label{section:integration}

Let $\Msc$ be holonomic $\Dsc_X$-module defined on $X := \C^n$.
In this section, we explain a method to
translate the computation of integrals
$\int^i_X \Msc$  to that of the restriction
$ H^i((\Dsc_X/(x_1 \Dsc_X + \cdots + x_n \Dsc_X) \otimes_{A_n}^L {\hat \Msc})$,
where $\hat \Msc$ is the Fourier transform of $\Msc$.  
Our discussion together with the algorithm of computing
the restriction in \cite{OakuAdvance2} prove 
the correctness of 
Procedure \ref{procedure:integral}
and consequently the correctness of steps 2, 3 and 4 of
Algorithm \ref{algorithm:cohomology-top}.

Let us denote by $\Omega_X^{i}$ the sheaf of regular (algebraic)
$i$-forms on $X$.  We use the notation 
$\partial = (\partial_1,\dots,\partial_n)$ with 
$\partial_i := \partial/\partial x_i$.  
Let us denote by ${\rm DR}(\Msc)$ the complex
$$
0 \longrightarrow \Omega_X^{0}\otimes_{\Osc_X}\Msc 
\stackrel{d}{\longrightarrow} \Omega_X^{1}\otimes_{\Osc_X}\Msc 
\stackrel{d}{\longrightarrow}
\cdots
\stackrel{d}{\longrightarrow} \Omega_X^{n}\otimes_{\Osc_X}\Msc 
\longrightarrow 0, 
$$
where $d$ is defined by
$$
d(dx_{k_1}\wedge \cdots \wedge dx_{k_i}\otimes u)
= \sum_{j=1}^n dx_j\wedge dx_{k_1}\wedge \cdots \wedge dx_{k_i}
\otimes (\partial_j u)
$$
for $u \in \Msc$.  Here we regard $\Omega^{i}\otimes_{\Osc_X}\Msc$ as
being placed at degree $i-n$.  
In particular, the cohomology groups of ${\rm DR}(\Dsc_X)$ are given by
$$
\Hsc^i({\rm DR}(\Dsc_X)) = \left\{\begin{array}{ll}
\Omega_X & \mbox{if $i=0$,}\\
0        & \mbox{otherwise.}\end{array}\right.
$$
Hence we have an isomorphism 
$$ 
\Omega_X {\otimes}_{\Dsc_X}^{L} \Msc 
\simeq {\rm DR}(\Dsc_X){\otimes}_{\Dsc_X}\Msc 
= {\rm DR}(\Msc).
$$
Since $\Omega_X^{i}\otimes_{\Osc_X}\Msc$ is a quasi-coherent $\Osc_X$-module
and $X$ is affine, we have 
$$
H^k(X,\Omega_X^{i}\otimes_{\Osc_X}\Msc) = 0 \quad (k > 0).
$$
Hence by using the standard argument for the sheaf cohomology, 
the integral is explicitly represented by a complex
$\int_X \Msc = R \Gamma(X; {\rm DR}(\Msc))$, which is equivalent to 
\begin{equation} \label{equation:deRham1}
0 \longrightarrow (\stackrel{0}{\wedge}\Z^n)\otimes_\Z M 
\stackrel{d}{\longrightarrow} (\stackrel{1}{\wedge}\Z^n) \otimes_\Z M
\stackrel{d}{\longrightarrow}
\cdots
\stackrel{d}{\longrightarrow}  (\stackrel{n}{\wedge}\Z^n) \otimes_\Z M
\longrightarrow 0, 
\end{equation}
where $M := \Gamma(X,\Msc)$ and 
$$
d(e_{i_1}\wedge\dots\wedge e_{i_k}\otimes u) = 
\sum_{j=1}^n e_j\wedge e_{i_1}\wedge\dots\wedge e_{i_k}\otimes (\partial_ju)
$$
with the unit vectors $e_1,\dots,e_n$ of $\Z^n$.

The Weyl algebra $A_n$ has a ring automorphism $\Phi$ defined by
$$ 
\Phi(x_i) =  -\partial_i,\quad
\Phi(\partial_i) =  x_i \qquad (i=1,\dots,n),
$$
This $\Phi$ naturally defines a new left $A_n$-module ${\hat M}:=\Phi(M)$, 
which is called the Fourier transform of $M$.  
Since $\Msc$ is holonomic, $M$ belongs to the Bernstein class 
of $A_n$-modules (cf.\ \cite[p.125]{BjorkBook}).  
Since the Bernstein class is invariant under the Fourier transform, 
we know that $\hat\Msc := \Dsc_X\otimes_{A_n}\Phi(M)$ 
is a holonomic $\Dsc_X$-module on $X$.   
By applying $\Phi$ to the complex (5.1), we obtain another complex
\begin{equation} \label{equation:deRham2}
0 \longrightarrow (\stackrel{0}{\wedge}\Z^n)\otimes_\Z \Phi(M) 
\stackrel{\delta}{\longrightarrow} (\stackrel{1}{\wedge}\Z^n) \otimes_\Z\Phi(M)
\stackrel{\delta}{\longrightarrow}
\cdots
\stackrel{\delta}{\longrightarrow}  (\stackrel{n}{\wedge}\Z^n) 
\otimes_\Z\Phi(M)
\longrightarrow 0, 
\end{equation}
where
$$
\delta(e_{i_1}\wedge\dots\wedge e_{i_k}\otimes u) = 
\sum_{j=1}^n e_j\wedge e_{i_1}\wedge\dots\wedge e_{i_k}\otimes (x_ju).
$$
Since the complexes (\ref{equation:deRham1}) and 
(\ref{equation:deRham2}) are isomorphic, we have only to
compute the cohomology groups of (\ref{equation:deRham2}).  
Here note that (\ref{equation:deRham2}) 
is a complex defining the restriction of $\hat\Msc$ 
to the origin of $X$.  
Thus, we have the following proposition.

\begin{proposition}
\label{proposition:6.1}
We have for any $i$, 
$$
H^i(X,DR(\Msc))\simeq H^i((A_n/(x_1A_n +\dots + x_nA_n))\otimes^L_{A_n}\hat M).
$$
\end{proposition}

Note that $\hat\Msc$ is specializable to the origin (i.e., a nonzero
$b$-function exists) since $\hat\Msc$ is holonomic (cf.\ \cite{KKholonomic3}).
Hence the cohomology groups of (\ref{equation:deRham2}) 
are computable by steps 2 -- 6 of Procedure \ref{procedure:integral}
as shown in \cite{OakuAdvance2}.
Thus each $\int_X^i\Msc$ is computable as a finite dimensional vector space
and we obtain the following proposition.

\begin{proposition}
\label{proposition:6.2}
Procedure \ref{procedure:integral} is correct.
\end{proposition}

The heart of Procedure \ref{procedure:integral} is the truncation
of the resolution with respect to a filtration defined by the
weight vector $w$
by a root of $b$-function \cite{OakuAdvance2}.
Let us briefly explain the idea.
Let $\Msc$ be a holonomic $\Dsc$-module and $b(s;x)$ be the $b$-function
(or indicial polynomial) along $x_1=0$.
Then,
$ x_1 \cdot \,: gr_{k+1}(\Msc)_p \rightarrow gr_k(\Msc)_p$
is bijective if $b(k;p)\not=0$.
Here $gr(\Msc)$ is the graded module associated with
the weight vector $w = (-1, 0, \ldots, 0; 1, 0, \ldots, 0)$.
Hence, in order to obtain the kernel and the image of the map $x_1$,
we may truncate the high degree part and the low degree part of the
filtration of $\Msc$ with respect to the weight vector $w$.
In order to obtain all the cohomology groups of the restriction, we need
a diagram chase to determine the degree of the truncation.
As to details, see \cite[section 5]{OakuAdvance} and \cite{OakuAdvance2}.

\bigbreak

By Proposition \ref{proposition:3.10}, Proposition \ref{proposition:4.2},
Theorem \ref{theorem:comparison} and Proposition \ref{proposition:6.2},
we obtain the following theorem and complete our proof
of Theorems \ref{theorem:constantSheaf} and \ref{theorem:locallyConstantSheaf}.

\begin{theorem}
Algorithms \ref{algorithm:cohomology-top} and \ref{algorithm:cohomology-top2}
are correct.
\end{theorem}

We close this section with the following theorem,
which generalizes Theorem \ref{theorem:locallyConstantSheaf}
under the condition that the coefficient sheaf is expressed 
in terms of $A_n$-module $M$.
Note that we did not require this condition in Theorem
\ref{theorem:locallyConstantSheaf}.

\begin{theorem}
Let $M$ be an $A_n$-module 
$ (A_n)^p/ I $ where $I$ is a left submodule of $(A_n)^p$.
We assume that $\Msc = \Dsc_X M$ is regular holonomic on $X$ and
also regular along ${\bf P}^n \setminus {\bf C}^n$ and that
the singular locus of $\Msc$ on $X$ is given by $f=0$
with a polynomial $f$.
The cohomology group
$ H^k(U, \Hom_{\Dsc_U}(\hol_U, \Dsc^{\rm an}_U\otimes_{A_n}M))$
is computable where
$U =\C^n \setminus V(f)$.
\end{theorem}

\begin{proof}
An algorithm to compute
$ M[1/f] $ is given in \cite[Tensor product and localization]{OakuAdvance2}
(see also \cite[Propositions 7.1 and 7.5]{OakuAdvance}).
The rest of the proof is same as that of 
Theorems \ref{theorem:constantSheaf} and
\ref{theorem:locallyConstantSheaf}
because $\Msc_{|_{f \not= 0}}$
can be regarded as a regular connection
(\cite[Theorem 2.3.2]{KKholonomic3}).
\end{proof}

\begin{example} \rm
If $\Msc$ is not regular holonomic on ${\bf P}^n$, then the comparison theorem
no longer holds.
For example, put
$$ \Msc = \Dsc_X/ \langle \pd{x}+ 2 x \rangle, \quad
   X = U = \C.$$
The operator
$\pd{x} + 2 x$ is not regular at $x=\infty$.
We can see that
$$H^0(X, DR(\Msc)) = \C, \ 
  H^{-1}(X, DR(\Msc)) = 0
$$
by applying our integration algorithm.
Now, take $\varphi \in \Hom_{\Dsc_X}(\Osc_X,\Msc^{\rm an})$.
We may assume that $f = \varphi(1)$ belongs to $\Osc^{\rm an}$ since
$\pd{x} = -2x$ in $\Msc^{\rm an}$.
We have  $\pd{x} f  = 0$ in $\Msc^{\rm an}$, which means that
$ f \pd{x} + f' \in \langle \pd{x} + 2 x \rangle$.
Then, we have $f'/f = 2x$ and hence 
$ f = \varphi(1) = c e^{x^2} \in \Msc^{\rm an}$
for a constant $c$.
Therefore, we have
$ \Hom_{\Dsc_X}(\Osc_X,\Msc^{\rm an}) \simeq \C$.
On the other hand,
$$ H^1(X,\C) = 0 \not= H^0(X,DR(\Msc)) \ \mbox{ and }
   H^0(X,\C) = \C \not= H^{-1}(X,DR(\Msc)).
$$  
\end{example}

\section{Computation of cohomology groups on the complement of an 
algebraic set when its algebraic local cohomology group vanishes
except for one degree}
\setcounter{equation}{0}
\label{section:acyclic}

The purpose of this section is to establish a connection between 
the de Rham cohomology of $\C^n$ with an algebraic set removed, 
and the integration of modules over the Weyl algebra.  
We use the algebraic local cohomology groups lying in between these two
objects.
The contents of this section except the last theorem
should be well-known to specialists.

Put $X$ be an $n$-dimensional non-singular algebraic variety over $\C$
and let $Y$ be an arbitrary algebraic set of $X$. 
The algebraic local cohomology group $\Hsc^i_{[Y]}(\Osc_X)$ 
(in the sense of Grothendieck) is the $i$-th derived functor of 
the functor $\Gamma_{[Y]}$ of taking the support (in the algebraic sense). 
This is a holonomic $\Dsc_X$-module (Theorem 1.4 of \cite{K1}). 
When $X = \C^n$, 
algorithms for computing the algebraic local cohomology groups 
have been given in \cite{OakuAdvance} 
for the case where $Y$ is of codimension one and \cite{Walther}, \cite{OakuAdvance2}
for the general case.  

\begin{proposition}
\label{proposition:7.1}
Let $\C_X$ be the constant sheaf on $X$ with stalk $\C$.
Then there is an isomorphism
$$
R\Gamma_Y(X,\C_X) \simeq \int_X R\Gamma_{[Y]}(\Osc_X)[-n],
$$
where $[-n]$ denotes the shift operator in the derived category. 
In particular, if $\Hsc_{[Y]}^i(\Osc_X)=0$ for $i \neq d$, then
for any $i\in\Z$, there is an isomorphism
$$
H_Y^i(X; \C_X) \simeq \int_X^{i-n-d}\Hsc_{[Y]}^d(\Osc_X).
$$ 
\end{proposition}

\begin{proof}
The algebraic and the analytic de Rham complexes are defined by
\begin{eqnarray*}
{\rm DR}(\Osc_X) &:=& \Omega_X {\otimes}_{\Dsc_X}^{L} \Osc_X,\\
{\rm DR}^{{\rm an}}(\Osc^{\rm an}_X) &:=& 
\Omega^{\rm an}_X {\otimes}_{\Dsc^{\rm an}_X}^{L} \Osc^{\rm an}_X,
\end{eqnarray*}
where $\Omega^{\rm an}_X$ denotes the sheaf of holomorphic $n$-forms.
Then we have
\begin{eqnarray*}
\int_X R\Gamma_{[Y]}(\Osc_X) 
&=& R\Gamma(X; \Omega_X {\otimes}_{\Dsc_X}^{L} R\Gamma_{[Y]}(\Osc_X))\\
&=& R\Gamma(X; R\Gamma_{[Y]}({\rm DR}(\Osc_X)))\\
&=& R\Gamma_Y(X; {\rm DR}(\Osc_X)).
\end{eqnarray*}
On the other hand, there are two distinguished triangles and a morphism
between them:
$$\begin{array}{ccccc}
R\Gamma_Y(X; {\rm DR}(\Osc_X)) &\longrightarrow &
R\Gamma(X; {\rm DR}(\Osc_X)) &
\stackrel{+1}{\longrightarrow} &
R\Gamma(X\setminus Y; {\rm DR}(\Osc_X))\\
\downarrow & & \downarrow & & \downarrow\\
R\Gamma_Y(X; {\rm DR}^{\rm an}(\Osc^{\rm an}_X)) &\longrightarrow &
R\Gamma(X; {\rm DR}^{\rm an}(\Osc^{\rm an}_X)) &
\stackrel{+1}{\longrightarrow} &
R\Gamma(X\setminus Y; {\rm DR}^{\rm an}(\Osc^{\rm an}_X)).
\end{array}$$
Here the vertical homomorphisms except the leftmost one are isomorphisms
by virtue of the comparison theorem of Grothendieck \cite{G}.  
Hence the leftmost vertical homomorphism is also an isomorphism.
Moreover the complex de Rham lemma implies
${\rm DR}^{\rm an}(\Osc^{\rm an}_X) = \C_X[n]$.
Consequently, we get 
$$
R\Gamma_Y(X; {\rm DR}(\Osc_X)) 
= R\Gamma_Y(X; {\rm DR}^{\rm an}(\Osc^{\rm an}_X)) 
= R\Gamma_Y(X; \C_X)[n].
$$
This completes the proof.
\end{proof}

  From the above proposition and the isomorphism 
$H^i(X\setminus Y; \C) \simeq H^{i+1}_Y(X; \C)$
(see, e.g., \cite[p.212, exercises 2.3]{HartshornBook}), we obtain 

\begin{corollary}
Assume $\Hsc_{[Y]}^i(\Osc_X)=0$ for $i \neq d$.
Then for any $i \geq 1$, we have an isomorphism
$$
H^i(X\setminus Y; \C_X) \simeq \int_X^{i-n-d+1}\Hsc^d_{[Y]}(\Osc_X).
$$
\end{corollary}

If $Y$ is non-singular, we can also relate the de Rham cohomology of 
$X\setminus Y$ to that of $Y$ itself:

\begin{corollary}
Assume that $Y$ is non-singular and of codimension $d$.  
Then, for any $i\geq 1$, there exists an isomorphism
$$
H^i(X\setminus Y; \C_X) \simeq H^{i+1-2d}(Y; \C_Y).
$$
Hence, $H^i(Y; \C_Y)$ is computable for any $i \geq 0$.
\end{corollary}

\begin{proof}
Let $\iota : Y \rightarrow X$ be the embedding.
Then by the Kashiwara equivalence 
(cf.\  \cite[p.34, Theorem 1.6.1]{HottaTanisaki} and \cite{K1}), 
we have an isomorphism
$\Hsc^d_{[Y]}(\Osc_X) = \iota_+\Osc_Y$.
Thus by using Proposition \ref{proposition:7.1}, we obtain
\begin{eqnarray*}
H_Y^i(X; \C_X) &\simeq& \int_X^{i-n-d}\Hsc^d_{[Y]}(\Osc_X)\\
&\simeq& \int_X^{i-n-d}\iota_+\Osc_Y\\
&\simeq& \int_Y^{i-n-d}\Osc_Y\\
&\simeq& H^{i-2d}(Y; \C_Y).
\end{eqnarray*}
Combining this with the preceding corollary, we are done.
\end{proof}

In \cite{Walther}, U.~Walther gave an algorithm to 
compute the local cohomology groups $\Hsc^k_{Y}(\Msc)$
with a Cech complex under the condition that
$\Msc$ is $(f_1 \cdots f_d)$-saturated.
Since $\Osc_X$ satisfies this condition,
we can compute
algebraic local cohomology groups $\Hsc_{[Y]}^i(\Osc_X)$ for any 
$i \geq 0$  where $Y := \{f_1 = \dots = f_d = 0\}$.
Another approach to compute algebraic local cohomology groups
of $\Msc$  with a resolution and 
without the condition of saturation is given in
\cite{OakuAdvance2}.
Thus, we have two algorithms for the next theorem.

\begin{theorem}
The cohomology groups $H^i(X\setminus Y; \C_X)$ for any $i\geq 0$ 
is computable
if $ Y = V(f_1, \ldots, f_d),\  f_i \in \Q[x_1, \ldots, x_n]$ 
and if
$\Hsc_{[Y]}^j(\Osc_X)$ vanishes except for one $j$.
\end{theorem}

This theorem generalizes Theorem \ref{theorem:constantSheaf} 
under the condition on vanishing of the local cohomology groups
$\Hsc_{[Y]}^j(\Osc_X)$. 
Note that if $d=1$, then this condition always holds.

\bigbreak
\bigbreak
\bigbreak

\def\at{\catcode`@=11{@}\catcode`@=12 } 
\bigbreak
\rightline{Toshinori Oaku}
\rightline{\tt oaku\at math.yokohama-cu.ac.jp} 
\rightline{Department of Mathematical Sciences, Yokohama City University}
\rightline{Seto 22-2, Kanazawa-ku, Yokohama, 236-0027, Japan}

\bigbreak
\rightline{Nobuki Takayama }
\rightline{\tt takayama\at math.kobe-u.ac.jp}
\rightline{Department of Mathematics, Kobe University}
\rightline{Rokko, Kobe, 657-8501, Japan.}

\end{document}